%% file: ffinal.tex
\documentclass[11pt,a4paper]{article}
\usepackage{a4wide}
\usepackage{amsmath,amssymb,amsthm,amsfonts}
\usepackage{psfig,epsfig,epic,eepic,graphicx}
\newtheorem{theo}{Theorem}
\newtheorem{lem}{Lemma}
\newtheorem{cor}{Corollary}
\newtheorem{prop}{Proposition}
\newtheorem{defi}{Definition}
\newtheorem{rem}{Remark}

\def\Im{{\cal I}m \,}
\def\pa{\partial}
\def\na{\nabla}
\def\curl{\hbox{curl }}

\def\dis{\displaystyle}
\def\mspace{{\medskip \noindent}}
\def\sspace{{\smallskip \noindent}}

\newcommand{\N}{{\mathbb N}}
\newcommand{\Z}{{\mathbb Z}}
\newcommand{\Q}{{\mathbb Q}}
\newcommand{\R}{{\mathbb R}}
\newcommand{\C}{{\mathbb C}}
\newcommand{\T}{{\mathbb T}}

\newcommand{\E}{{\mathbf e}}
\newcommand{\eps}{{\varepsilon}}

\title{Remarks on Boundary Layer Expansions }
\author{David {\sc Gerard-Varet}, Thierry {\sc Paul}
\footnote{D.M.A, U.M.R. 8553, E.N.S., 45
   rue d'Ulm, 75005 Paris}}
\date{}
\begin{document}
\maketitle
\tableofcontents

\begin{abstract}
A systematic mathematical methodology for derivation of boundary layer
expansions is presented. An explicit calculation of boundary layer sizes is
given and proved to be coordinates system independent. It relies on
asymptotic properties of symbols of operators. Several examples, including
the quasigeostrophic model, are discussed. 
\end{abstract}

\section{Introduction and main result} \label{intro}
Boundary layers appear in various physical contexts, such as fluid
mechanics, thermodynamics, or ferromagnetic media. From the mathematical
point of view, they are related to singular perturbation problems in
bounded domains.  The singular perturbation is due to the presence of small 
parameters  in the dimensionless governing equations. For instance, in
magnetohydrodynamics (MHD),  
the so-called Rossby number (ratio between the  angular
velocity of the fluid and the angular velocity of the Earth) can be as
low as $10^{-7}$. Similarly, the Prandtl number (ratio between
hydrodynamic and magnetic diffusivity) is about $10^{-6}$ in the Earth's core.
In the interior of the domain, these small parameters lead to some reduced
dynamics of the equations. For instance, in highly rotating fluids, the
velocity field does not vary along the rotation axis, see textbook
\cite{Gree}. This reduced dynamics is often incompatible with boundary
conditions. This yields boundary layers, 
in which the solutions of the equations have strong gradients, in order to
satisfy the boundary conditions (typically a Dirichlet condition at a rigid
surface). 

\medskip
In a  formal setting, boundary layer problems are connected  to
systems of the type
\begin{equation} \label{general}
{\cal A}^{\eps}(x, D_x) u^{\eps} \: + \: {\cal Q}^{\eps}\left(u^{\eps} \right) 
 \: = \: f^{\eps}, \: \: 
 x \in \Omega,
\end{equation}
where $\eps \in \R^d$ describes the small parameters of the system, $\Omega
\subset \R^n$ is the domain, 
$u^{\eps}(x)$  the unknown, $f^{\eps}(x)$ the data, 
${\cal A}^{\eps}$ is a linear differential operator (most often of elliptic
or parabolic type),  and ${\cal Q}^{\eps}$ the nonlinear part. One must of
 course supply these equations with appropriate  boundary conditions. 

\medskip
The basic idea, which underlies the mathematical study of all boundary
layers, is that the solution $u^{\eps}$ of \eqref{general} should satisfy an
asymptotics of the form  
\begin{equation} \label{asymptotique}
\begin{aligned}
& u^{\eps}(x) \:  \sim \: u(x) \: + \: u_1^{bl}\left(x, \frac{d(x, \partial
\Omega)}{\alpha_1(\eps)}\right) \: + \ldots \: + \: u_r^{bl}\left(x, \frac{d(x, \partial
\Omega)}{\alpha_r(\eps)}\right), \quad \eps \rightarrow 0, \\
&  u^{bl}_i\left(\theta_i\right) \:  \rightarrow 0, \quad
\theta_i \rightarrow +\infty, \quad \forall i .
\end{aligned}
\end{equation}
This means that $u^{\eps}$ should have a regular part, depending on $x$,
but also a singular part, depending on stretched  variables $\theta_i =
d(x, \partial \Omega)/\alpha_i(\eps)$, with $\alpha_i(\eps)
\rightarrow 0$ as $\eps \rightarrow 0$. This singular part should be
localized near the boundary and express the strong gradients
of  boundary layers. Broadly speaking, the aim  of boundary
layers studies is to answer the following questions:
\begin{description}
\item[i)] What are the possible $\alpha_i(\eps)$ (the boundary layer sizes)?
\item[ii)] What are the possible profiles $u$, $u^{bl}_i$ ?
\item[iii)] Is the asymptotics \eqref{asymptotique} correct (in a sense to be
determined) ?
\end{description}
The first two questions are related to  the derivation of boundary layers,
whereas the third one is connected to stability issues.

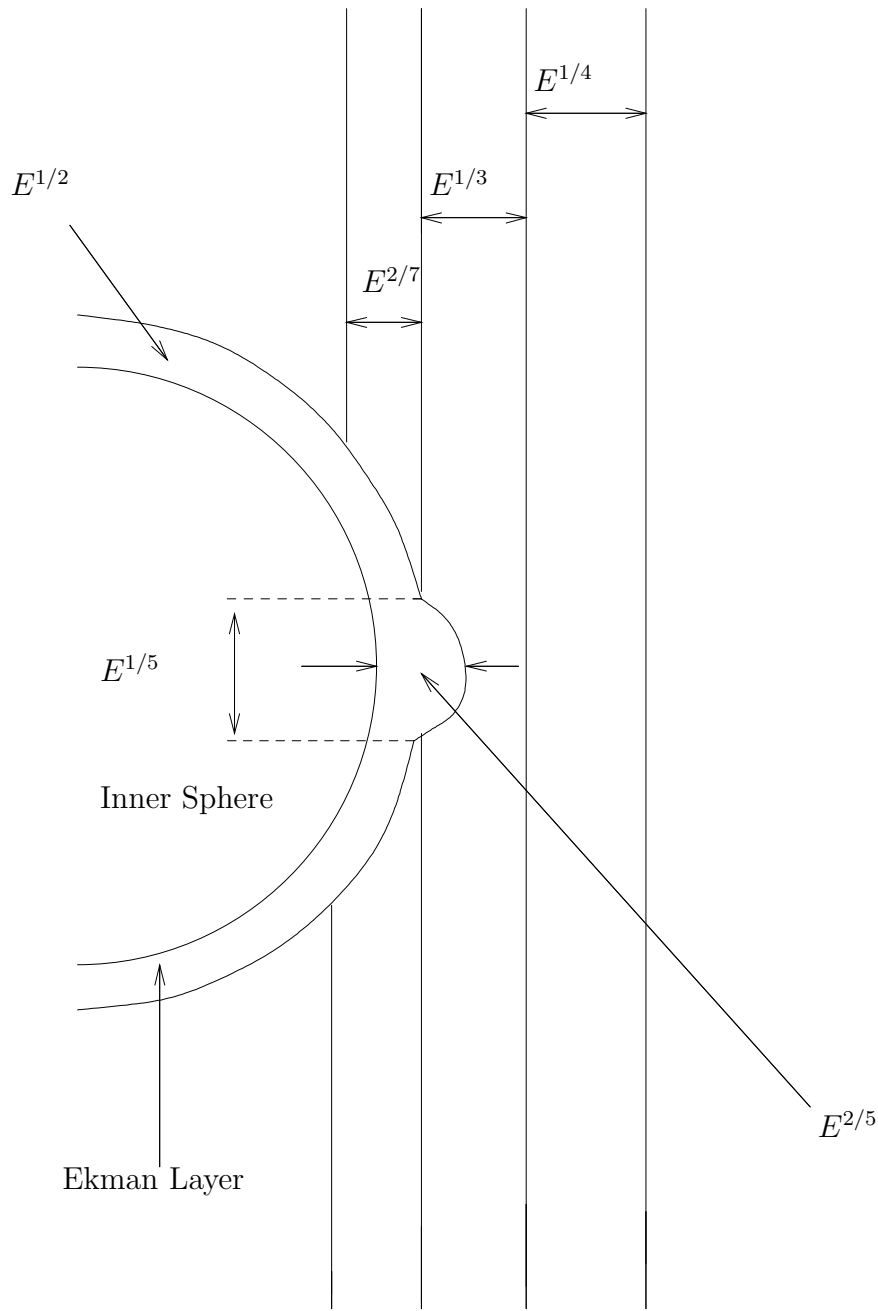
\begin{figure} \label{figure}
\begin{center}
\label{stewartson}
\input{sphere.eepic}
\caption{Boundary layers of rotating fluids, near a sphere and 
at the circumscribing cylinder (following Stewartson \cite{Ste2})}
\end{center}
\end{figure}

\medskip
The aim of this note is to give some insight into the derivation problem. In
some cases, the derivation is quite easy, and the stability analysis is the
most difficult part (see for instance \cite{GiSe, Mezu} 
on viscous perturbations of
hyperbolic systems). But in most situations of physical interest, it may
involve a variety of length scales and equations. A typical example is the
description of water in a highly rotating tank. 
If the water is at rest
in the rotating frame attached to the tank, it is well modeled by Stokes
equations with Coriolis term:
\begin{equation} \label{FT}
\begin{aligned}
\E \times u \: + \: \nabla p \: - \: E \Delta u \: & =  \: f, \\ 
  \nabla \cdot u \:  & =  \: 0,
\end{aligned}
\end{equation}
 where $\E = (0,0,1)$ is the rotation vector, and $E$ is a small parameter
 called the Ekman number. For such system (with appropriate 
source term $f$ and boundary conditions),  
the structure of the solutions near  a flat
horizontal boundary is  understood: there is  a 
boundary layer of  size $E^{1/2}$, the Ekman layer (see \cite{Ekm}).
 In this case, one can even carry  an  analysis of the full 
Navier-Stokes equations with rotation: we refer to \cite{GrMa, Mas1,
 ChDeGaGr3}  among others.   
But in more elaborate geometries, following the articles
 of Stewartson \cite{Ste1, Ste2}, 
many other layers develop: for instance, between
two concentric spheres with slightly different rotation speeds, boundary
layers of size   $E^{1/2}$, 
$E^{1/3}$,  $E^{1/4}$, $E^{2/5}$, \ldots are expected near the inner sphere
 and the cylinder circumscribing it (see figure \ref{figure}). 

\medskip
Note that equations \eqref{FT} are linear. Generally, the
derivation of the boundary layers does not involve the nonlinearity ${\cal
Q}^{\eps}$ of \eqref{general}: it  only matters in stability
questions.
 In short, one can say that the nonlinear term does not create the
boundary layer, but may destabilize it. Therefore, we 
restrict ourselves to linear equations 
\begin{equation} \label{linear}
 {\cal A}^{\eps}(x, D_x) \, u^{\eps} \: = \: f^{\eps}. 
\end{equation}

\medskip
Classically, the problem of the derivation is tackled through one of the
following ways:
\begin{description}
\item[i)] an explicit calculation, where the exact solution is  computed and 
 expanded. However, such
technique is restrictive (an analytical computation is rarely tractable) 
and often tedious (see for instance \cite{Ste2}, with Bessel functions).   
\item[ii)] the so-called ``method of matched asymptotics'', see  
\cite{Eck, Vand, Zey} for a general presentation. 
The idea of this method is to patch
two asymptotic expansions: an ``outer'' regular one, far from the boundary,
and an ``inner'' singular one, close to it. Such expansions must coincide
in an overlap domain, which  provides
 boundary conditions for both outer and inner profiles. The structure of
the inner expansion (including the layer sizes) is sought by trial, through
the principle of least degeneracy ({\it c.f.} \cite{Eck} 
for a detailed explanation). This method has been applied with success to 
various physical systems, including singularities near edges.  
However, the determination of layer sizes leads to heavy
computations, and must often be supplied with refined physical arguments
 (see \cite{Vande} on rotating fluids). 
\end{description}
A new  method of derivation  has been recently 
 introduced in \cite{Ger4} (see also the proceedings  \cite{GeGr}). 
It is based on
 the asymptotic analysis of the symbol  ${\cal A}^{\eps}(x, \xi)$ 
 as $\eps$ goes to zero, for  $x \in \partial \Omega$. 
Briefly, up to use local coordinates near $x \in \partial
 \Omega$,  one can consider  
that equation \eqref{linear} holds in a neighborhood of $x = (x', 0)$
 in $\R^n_+$. Then, the leading idea of the derivation is to carry a 
Fourier-Laplace analysis: 
 boundary layer sizes near $x$ are  deduced from  modal solutions 
\begin{equation} \label{modes}
u^{\eps}(\tilde{x}) \: = \: \exp\left(i \xi^{\eps} \cdot \tilde{x}\right)
 \, V^{\eps}, \quad V^{\eps}
\neq 0, 
\end{equation}
of the equation with frozen coefficients
\begin{equation*} 
{\cal A}^{\eps}(x,\partial_{\tilde{x}}) \,  u^{\eps} \: = \: 0.
\end{equation*}
In other words, one must consider the characteristic manifold of ${\cal
A}^{\eps}$,
$$ \sigma\left(a^{\eps}\right) = \left\{ \left(x,\xi^{\eps}\right),
 \: a^{\eps}(x, \xi^{\eps}) = 0 \right\}
$$   
where $a^{\eps}(x,\xi) = \det {\cal A}^{\eps}(x, \xi)$. 
Broadly speaking, if $\displaystyle  
\xi = (\zeta, \xi_n)$ is the dual variable of 
$\displaystyle x = (x', x_n) \in \R^n_+$,
boundary layers  correspond to $\left(x,\xi^{\eps}\right)$
 in $\sigma\left(a^{\eps}\right)$
 that satisfy
$$ \Im \xi_n^{\eps} \: \rightarrow +\infty, \quad \eps \rightarrow 0. $$
The size of the boundary layers is then  given by $\displaystyle 
\left| \Im \xi_n^{\eps} \right| ^{-1}$. In most cases, boundary layer
equations follow, by appropriate rescaling of the symbol. 
This method has been used in \cite{Ger4} at a fully formal level. It
has been applied with efficiency to various  geophysical systems, 
including rotating fluids or MHD. It has allowed to recover the main
boundary layers   of the physical literature, with very few and simple
algebraic computations (see \cite{Ger4} for all details).

\medskip
The aim of the present paper is to provide this formal method with some
 mathematical basis. We will  limit ourselves to  scalar  equations. 
In this reduced (but still large!) setting, we will show how the microlocal  
analysis of \eqref{linear} is linked to  boundary layer expansions
\eqref{asymptotique}.   We will 
 present a set of  conditions (assumptions (H1) to (H5)) that
 ensures the existence of such expansions, {\em 
with clearly identified boundary
 layer sizes and equations}. Hypothesis (H1) to (H5)
 involve symbols derived from  $a^{\eps}$. They will be shown to be
 {\em intrinsic}, although these symbols depend on the choice of local
 coordinates.

\medskip
Let us  specify the framework of the study. Let $X$ be a smooth
Riemannian manifold,  $n=\dim X$. In all situations of physical interest, 
$X$ will be of the type $\T^{n_1} \times \R^{n_2}$, $n_1 + n_2 = n$, where
the  torus $\T^{n_1}$ models periodic boundary conditions in $n_1$ physical
variables. The Riemannian metric will  be the Euclidean metric
induced by $\R^{2 n_1} \times \R^{n_2}$. Let $\Omega$ a domain
of $X$. We denote $Y=\pa \Omega$ the boundary
of $\Omega$, which is a $n-1$ dimensional submanifold of $X$. We assume
that $Y$ has a finite number of connected components, $Y_1$ to $Y_N$. We  endow
$Y$ with the Riemannian structure induced by $X$.
We assume that there exist a tubular
neighborhood of $Y$ in $X$,  $T > 0$, and  a
smooth diffeomorphism
$$  \: {\cal T} \approx Y \times (-T, T), \quad x \approx (y,t), $$
such that $t > 0$ in ${\cal T} \cap \Omega$, $\: 
t < 0$ in ${\cal T} \cap (X \setminus \bar{\Omega})$. 
For instance, such assumption is satisfied when $X = \R^n$ and $\Omega$ is a
bounded open subset. 
{\em In the whole sequel, we identify ${\cal T}$  and $Y \times (-T, T)$}.

We consider the  scalar problem,
\begin{align}
 \label{elliptique}
&  a^{\eps}(x, D_x) \, u^{\eps} \: = \:  f, \quad x \in \Omega, \\
 \label{elliptiquebc}
& b^\eps_{l}(x, D_x) \,  u^\eps \: = \: g_l,
\quad x \in Y,  \quad l \in \{ 1, \dots, L \},
\end{align}
where $\dis a^{\eps}(x, D_x)$, respectively $\dis b^\eps_{l}(x, D_x)$,
  is a dif\-fe\-ren\-tial operator depending
smooth\-ly on  
$\dis x \in \Omega \cup {\cal T}$, respectively $x \in Y$,  and poly\-no\-mial
 in $\eps$. 
The source term $f$, respectively $g_l$, is smooth on $\dis 
\Omega \cup {\cal T}$, resp. on $Y$. Let $m = \deg(a^\eps) > 0$. 

In the region $Y \times [0, T')$, $T' > 0$ small enough,  
we can write
$$ a^\eps(x, D_x) \: =  \: a^\eps(x, D_y, D_t) \: =  \: 
 \sum_{j=0}^m \, a^{\eps}_j(x, D_y) \, D_t^j $$ 
where $a^{\eps}_j$ is a differential operator on $Y$, of order $m-j$, 
with coefficients smooth in $x$, polynomial in $\eps$. 
We assume that the leading coefficient reads
\begin{equation*} 
a_m^{\eps}(x) \: = \: \eps^M \, a_m(x), \quad M \in \N, \quad a_m(x) \neq
0, \: \forall x. 
\end{equation*}
 
\medskip
 Before we state our main result, we still need to precise what we mean by 
boundary layer expansion.  Let  $\varphi$ be a cut-off near the boundary.
 Precisely, $\displaystyle \varphi \in {\cal
C}^\infty(\overline{\Omega})$, $\displaystyle 
\varphi = 1$ on $Y \times [0, T/4)$, and
$\displaystyle 
\varphi = 0$ on $\overline{\Omega} \setminus \left( Y \times [0, T/2)\right)$. 
  
\begin{defi}
A family of functions $u^\epsilon(x)$  is {\em of boundary layer type }
if it reads
 \begin{equation} \label{ansatz1}
u^\eps(x) \: = \: u^\eps_r(x) \: + \: \varphi(x) v^\eps_{bl}(x) 
\end{equation}
where the regular part $ u^\eps_r$ and the singular part $v^\eps_{bl}$ have
the following asymptotic expansions, for some positive $\delta$ and
$\gamma^i_j$:
\begin{equation} \label{ansatz2}
u^\eps_r(x) \: \sim \: \sum_{k=0}^{\infty} \eps^{\delta k} \, u^k(x), \quad 
 u^k \in {\cal C}^\infty\left(\overline{\Omega}\right),
\end{equation}
uniformly on every compact subset of $\overline{\Omega}$, and 
 and 
\begin{equation} \label{ansatz3}
 v^\eps_{bl}(y,t) \: \sim \: \sum_{k=0}^{\infty} \eps^{\delta k}
\sum_{\substack{1 \le i \le N \\ 1 \le j \le r_i}} \mathbf{1}_{Y^i}(y) \,
v^{i,k}_j\left(y, \frac{t}{\eps^{\gamma^i_j}}\right), 
\end{equation}
 uniformly on
every compact subset of $Y \times [0, T)$, 
where $v^{i,k}_j = v^{i,k}_j(y, \theta)$ are not all zero and belong
 to ${\cal X}^i = 
{\cal S}\left( \R^+_t ; \, {\cal C}^\infty(Y^i)\right)$. 
\end{defi}
We recall that ${\cal S}$ stands for the Schwartz space of fast decreasing
functions. Spaces $\chi^i$ express the localization of boundary layer
profiles. Note that in our definition, the singular part is non zero (at
least one $v^{i,k}_j$ is non zero). Our main result resumes to    
\begin{theo}\label{theo1}
Suppose that assumptions (H1) to (H5), {\it c.f.}  next sections, hold.
  Then,  the system  
  \begin{equation}
  \left\{
  \begin{aligned} 
& a^{\eps}(x, D_x) \, u^{\eps} \: = \:  f, \quad x \in \Omega, \\
& b^\eps_{l}(x, D_x) \,  u^\eps \: = \: g_l,
\quad x \in Y,  \quad l \in \{ 1, \dots, L \},
  \end{aligned}
  \right.
  \end{equation}
  has a solution of boundary layer type, modulo $\eps^\infty$
  uniformly on every compact subset.
\end{theo}
Theorem \ref{theo1} is the rigorous translation of the formal method
described in article \cite{Ger4}. In the following, we will present assumptions
(H1) to (H5), which allow to build solutions of boundary layer
type. Broadly, assumptions (H1) to (H3) (see section
\ref{singularperturbation}) relate to the determination of exponents
$\gamma^i_j$. Assumption (H4) (section \ref{blsizes}) is linked to the
derivation of boundary layer profiles $v^{i,k}_j$. Finally, (H5) (section
\ref{blexpansions}) is connected to the construction of the regular part
$u^\eps_r$. The proof of theorem \ref{theo1} follows. Application to the
quasigeostrophic equation is given in section \ref{sectionqg}.

\section{The singular perturbation} \label{singularperturbation} 
\subsection{Local characteristic manifolds and assumptions (H1)-(H2)}
The aim of the paper is to describe the structure of the solutions
 $u^{\eps}$ of \eqref{elliptique}  near the boundary.
 Let $y \in Y$, and $\displaystyle \left({\cal O}, 
\chi=(x', x_n) \right)$ a local chart in X around $y$, with 
$$ x_n(y) = 0, \quad \frac{\pa x_n}{\pa t}(y) \neq  0. $$
The local coordinates $(x', x_n)$ lie in $\R^{n-1} \times \R$. We denote by 
$a^{\eps}_{\chi} = a^{\eps}_\chi(x, \xi)$ the symbol of the operator 
$a^{\eps}$ in this local chart: precisely,
$$ a^{\eps}_{\chi}(x, \xi) \: := \: e^{-i\chi(x) \cdot \xi} \left[
 a^{\eps}(x, D_{\tilde{x}}) e^{i \chi(\tilde{x})\cdot  \xi} 
\right]_{\vert \, \tilde{x} = x}. $$
Our idea is to deduce the singular structure of $u^{\eps}$ from the
local symbols $a^{\eps}_\chi$, precisely from their characteristic
manifolds. We are interested in modal solutions that are singular with
respect to $\eps$ in the direction normal to the boundary. This means we 
wish to consider wavevectors $\xi^{\eps} = (\zeta, \xi^{\eps}_n)$ 
 such that 
$$ | \xi_n^{\eps}  | \rightarrow +\infty, \quad \eps \rightarrow 0. $$
We have the following 
\begin{prop} \label{holom1}
For all $\: \zeta \in \R^{n-1}$, 
\begin{description}
\item[i)] the m roots of $\displaystyle 
a^\eps_\chi(y, \zeta, \cdot)$ can be
written as m functions $\: \xi_1^{\zeta}(\eps)$, $\ldots$, 
$\xi_m^{\zeta}(\eps)$ for $\eps > 0$ small enough.
\item[ii)] There exists $p = p(\zeta) \in \N^*$ such that, for $i=1,
\ldots, m$, $\xi_{i}^{\zeta}(\eps^p)$  has an extension
 meromorphic in $\eps$. 
\end{description} 
\end{prop}  
This proposition follows from standard results of complex analysis. For
fixed $\zeta$,  $a^\eps_\chi(y, \zeta, \cdot)$ is a polynomial in
$\xi_n$  with coefficients holomorphic in $\eps$, 
and we refer to Kato \cite{Kat} for a detailed study. 
\begin{cor}
For all $\: \zeta \in \R^{n-1}$, $1 \le i \le m$, one of
the two following possibilities occurs:
\begin{enumerate} \label{holom2}
\item $\xi_{i}^{\zeta}(\eps)$ has a limit 
as $\eps > 0$ goes to zero, that we denote $\xi_{i,0}^{\zeta}$.
\item there exists a unique $p_i = p_i(\zeta) \in \Q^+_*$ such that
$\eps^{p_i} \, \xi_{i}^{\zeta}(\eps)$ has a non-zero limit  as $\eps >
0$ goes to zero, that we denote  $\eta_{i,0}^{\zeta}$ .
\end{enumerate}   
\end{cor} 
 Note that if the leading coefficient of  
$a^\eps_\chi(y, \zeta, \cdot)$ does not vanish 
when $\eps \rightarrow 0$, the second possibility does not occur. Roots
with  singular behaviour exist only if the highest order term goes
to zero with $\eps$, which is typical of a singular perturbation. 
 
\medskip
We denote by $s = s(\zeta)$ the number of roots 
$\displaystyle 
\xi_{i}^{\zeta}$ satisfying 2. Several indices $i$ may correspond to the
same value of $p_i$. Let $\displaystyle \gamma_1$, \ldots, $\gamma_r$,
 $r \le s$, be the distinct
values of $p_i$. For $1 \le j \le r$, we call $m_j$ the number of indexes
$i$ such that $\displaystyle p_i = \gamma_j$. We make the following assumption:

\mspace
{\bf (H1)   There exists $R = R(y) > 0$, such that 
 for all $|\zeta| > R$, the values of 
 $r$, $\: (\gamma_1, m_1)$ up to  $\: (\gamma_r,  m_r)$  are independent of 
$\zeta$.}    
  
\medskip
This also implies that $s = \sum  m_j$ is independent on $\zeta$. 
 Note that, using notation $\gamma_r$, we 
implicitly assume  $r \ge 1$. 
When $r = 0$,  the statements in the sequel 
become empty, so that we do not pay attention to this case.

\medskip
Assumption (H1) allows to state 
\begin{prop}  \label{zeros}
Let $ | \overline{\zeta} | > R$. There exists
$\displaystyle {\cal W} = {\cal W}\left(\overline{\zeta}\right)$
a neighborhood of $\overline{\zeta}$, and
constants $\displaystyle 
\eps_0 = \eps_0\left(\overline{\zeta}\right) >
0$,  $\displaystyle \delta = \delta\left(\overline{\zeta}\right) > 
0$, such that: for all $\zeta$ in $\displaystyle  {\cal W}$, for all
 $0 < \eps \le
\eps_0$, the roots of $a^{\eps}_{\chi}(y, \zeta, \cdot)$ can be divided
into $r+1$ disjoint sets $\displaystyle {\cal Z}_j^{\eps}(\zeta)$, 
$0 \le j \le r$
satisfying: 
\begin{description}
\item[i)] 
$\mbox{ card } {\cal Z}_0^{\eps}(\zeta) = m - s$, and for all $\xi_n$ in 
${\cal Z}_0^{\eps}(\zeta)$, 
$$ 0 \: \le \: | \xi_n | \: \le \: \delta^{-1}. $$
\item[ii)] 
For all $1 \le j \le r$, $\mbox{ card } {\cal Z}_j^{\eps}(\zeta) = m_j$,
and for all $\xi_n$ in  ${\cal Z}_j^{\eps}(\zeta)$,
$$ \frac{\delta}{\eps^{\gamma_j}} \: \le \: | \xi_n | \: \le \:
\frac{\delta^{-1}}{\eps^{\gamma_j}}.$$
\end{description}
\end{prop}
\begin{rem} \label{remark1}
Suppose proposition \ref{zeros} holds. Let $\zeta \in {\cal W}$.  By
propositions \ref{holom1} and \ref{holom2}, for $\eps$ small enough, 
the roots of $a^{\eps}_{\chi}(y, \zeta, \cdot)$ can be written 
$\xi_i^{\zeta}(\eps)$, $1 \le i \le m$. Up to reindex the roots, they
satisfy:
\begin{align*}
 \xi_i^{\zeta} \: & \xrightarrow[\eps \rightarrow 0]{} \: 
\xi_{i,0}^{\zeta}, \quad s+1 \le i \le m, \\
 \eps^{\gamma_j} \, \xi_i^{\zeta} \: & 
\xrightarrow[\eps \rightarrow 0]{} \: 
\eta_{i,0}^{\zeta} \neq 0, \quad 1 \le j \le r, \quad p_i = \gamma_j.
\end{align*}
Therefore, there exists $\eps_0^{\zeta}$ such
that , for all $0 < \eps < \eps_0^{\zeta}$, 
$$ {\cal Z}_0^{\eps}(\zeta) \: = \: \left\{ \xi_i^{\zeta}(\eps), \quad 
s+1 \le i \le m \right\}, $$
$$ {\cal Z}_j^{\eps}(\zeta) \: = \: \left\{ \xi_i^{\zeta}(\eps),
\quad p_i = \gamma_j \right\} $$
for $1 \le j \le r$. We emphasize that in proposition \ref{zeros}, the sets $
 {\cal Z}_j^{\eps}(\zeta)$ are defined for $0 < \eps \le \eps_0$, where
$\eps_0$ depends only on $\overline{\zeta}$.
 This is not the case for the 
$\xi_i^{\zeta}(\eps)$, defined only up to $\eps_0^{\zeta}$,
which  depends {\it a priori} on $\zeta$. Moreover, 
 up to consider a smaller $\eps_0$, sets ${\cal Z}^{\eps}_j(\zeta)$
satisfying i) and ii) will be automatically disjoint. 
\end{rem}
{\bf Proof:} Let $|\zeta| > R$. With the
notations  of  remark
\ref{remark1}, for all $0 < \eps \le \eps_0^{\zeta}$, we can factorize 
\begin{equation} \label{factor}
 a^{\eps}_{\chi}(y, \zeta, \xi_n) \: = \:  a^{\eps, \zeta}_0(\xi_n) \:
\prod_{j=1}^r \,
a^{\eps,\zeta}_{j}\left(\eps^{\gamma_j} \xi_n \right), 
\end{equation}
where 
\begin{align*}
 a^{\eps, \zeta}_0(\xi_n) & \: = \: \eps^M \, \eps^{-\sum \gamma_j
m_j}\,  C_{\chi} \, 
\prod_{i=s+1}^m \, \left( \xi_n \, - \, \xi_i^{\zeta}(\eps) \right), \\
 a^{\eps,\zeta}_{j}(\eta) & \: = \:\prod_{p_i = \gamma_j}  \, \left( 
 \eta \, - \, \eps^{\gamma_j} \xi_i^{\zeta}(\eps)  \right), \: \mbox{ 
for } 1 \le j \le r,
\end{align*}
and $\displaystyle  C_{\chi} \: := \: \left
( \frac{\pa x_n}{\pa t}(y) \right)^m \, a_m(y)$. 
We then define 
$$ a^{\eps}_{r, \chi}(\zeta, \eta) 
\: := \: \eps^{-\beta_r} \: a^{\eps}_{\chi}\left(\zeta, 
\frac{\eta}{\eps^{\gamma_r}}\right),  $$
with $\beta_r = M - m \gamma_r$. 
It is a polynomial in $\eta$, whose coefficients are functions of the type 
$$  \sum_{k = -k_0}^{k_0} b_k(\zeta) \, \eps^{\alpha k},$$
 for some fixed $k_0 \in \N$, $\alpha > 0$,  and  polynomials
symbols  $b_k$. 

\medskip
Using \eqref{factor}, it is easily seen that: for all $\zeta$, locally
uniformly in $\xi_n$, 
$$ a^{\eps}_{r, \chi}(\zeta, \eta) 
\: \xrightarrow[\eps \rightarrow
0]{} \: C_{\chi} \, 
\eta^{m-m_r} \, \prod_{p_i = \gamma_r} \, \left( \eta \, - \,
\eta_{i,0}^{\zeta} \right). $$
It is {\it a priori} pointwise convergence in $\zeta$. However, the
pointwise convergence of the coefficients, of the type 
$$   \sum_{k = -k_0}^{k_0} b_k(\zeta) \, \eps^{\alpha k}
 \:  \xrightarrow[\eps
\rightarrow 0]{} \: b(\zeta) $$ 
holds if and only if $b_k = 0$ for $k < 0$, and  $b_0 = b$.
This implies that the coefficients of 
$$ a^{0}_{r, \chi}(\zeta, \eta) \: := \:  C_{\chi} \: 
\eta^{m-m_r} \, \prod_{p_i = \gamma_r} \, \left( \eta \, - \,
\eta_{i,0}^{\zeta} \right) $$
are polynomials,
 and that the convergence is in fact locally uniform in $\zeta$. 

\medskip
Let $\displaystyle
 | \overline{\zeta} | > R$. Let
$\Gamma_r$  a curve enclosing the $\displaystyle
 \eta_{i,0}^{\overline{\zeta}}$  and no other root
of $\displaystyle a^0_{r, \chi}(\overline{\zeta}, 
\cdot)$. As this
 polynomial has smooth
coefficients, Rouch\'e's theorem yields a neighborhood ${\cal W}$ of $
\displaystyle 
\overline{\zeta}$ such that, 
for all $\zeta$ in ${\cal W}$, 
$\displaystyle \Gamma_r$ encloses the $\displaystyle 
\eta_{i,0}^{\zeta}$ (and no other
root of $a^0_{r, \chi}(\zeta, \cdot)$ ). 
Now, the convergence of
$\displaystyle a^{\eps}_{r, \chi}$ to $a^0_{r, \chi}$ is locally uniform in 
$\zeta$. Still by Rouch\'e's theorem, and up to take smaller
 ${\cal W}$, we get: for all $\zeta$ in
$\displaystyle {\cal W}$, $\Gamma_r$ encloses  exactly $m_r$ roots of
$a^{\eps}_{r,\chi}(\zeta, \cdot)$. 

\medskip
Back to $a^{\eps}_{\chi}$, this yields a $\delta =
 \delta(\overline{\zeta}) > 0$, such that $a^{\eps}_{\chi}(y, \zeta,
 \cdot)$ has 
$m_r$ roots satisfying 
$$ \frac{\delta}{\eps^{\gamma_r}} \: \le \: | \xi_n | \: \le \:
\frac{\delta^{-1}}{\eps^{\gamma_r}}. $$
This provides us with the set ${\cal Z}^{\eps}_r(\zeta)$. 

\medskip
The construction of the other sets is made inductively, using the
polynomials 
$$ a^{\eps}_{j, \chi}(\zeta, \eta) \: := \: \eps^{-\beta_j} 
\: a^{\eps}_{\chi}\left(y,
\zeta, \frac{\eta}{\eps^{\gamma_j}}\right),$$
with $\beta_j = M - m \gamma_j + \sum_{k \ge j+1} (\gamma_j - \gamma_k)$. 
 We only indicate how to build ${\cal Z}^{\eps}_{r-1}(\zeta)$. 
The general induction argument is left to the  reader. 

\medskip
The coefficients  of the  polynomial  $a^{\eps}_{r-1, \chi}(\zeta, \cdot)$ 
are of the same type as $a^{\eps}_{r, \chi}(\zeta, \cdot)$. Using again
\eqref{factor}, we have this time, locally uniformly in $\zeta, \eta$, 
\begin{align*}
& a^{\eps}_{r-1, \chi}(\zeta, \eta) \: \xrightarrow[\eps \rightarrow
0]{}  a^{0}_{r-1, \chi}(\zeta, \eta),  \\
& a^{0}_{r-1, \chi}(\zeta, \eta) \: := \: 
C_{\chi} \, \eta^{m-m_r-m_{r-1}} \, \prod_{p_i = \gamma_{r-1}} \, 
\left( \eta \, - \, \eta_{i,0}^{\zeta} \right) \, \left( 
 \prod_{p_i = \gamma_r}  -\eta_{i,0}^{\zeta} \right), 
\end{align*}
and $a^{0}_{r-1, \chi}$ has smooth symbols in $\zeta$ as coefficients. 
Note that $$  \prod_{p_i = \gamma_r} \, -\eta_{i,0}^{\zeta} \:
\neq \: 0 $$
 by definition of the $\eta_{i,0}^{\zeta}$. 
Reasoning as above with Rouch\'e's theorem, up to reduce ${\cal W}$,
$\eps_0$ and  $\delta$, we find a curve $\Gamma_{r-1}$
 (independent on $\zeta$ in
${\cal W}$, $\eps \le \eps_0$) which encloses exactly 
$m_{r-1}$ roots of $a^{\eps}_{\chi}(y, \zeta, \cdot)$, all satisfying 
$$ \frac{\delta}{\eps^{\gamma_{r-1}}} \: \le \: | \xi_n | \: \le \:
\frac{\delta^{-1}}{\eps^{\gamma_{r-1}}}, \quad \zeta \in {\cal W}. $$
They define the set ${\cal Z}^{\eps}_{r-1}(\zeta)$. $\Box$

\medskip
Let $ \overline{\zeta}$, ${\cal W}$ and $\eps_0$
as in proposition \ref{zeros}. For all $\zeta$ in ${\cal W}$, for all
$0 < \eps \le \eps_0$, we can factorize the symbol as 
\begin{equation} \label{factor2}
 a^{\eps}_{\chi}(y, \zeta, \xi_n) \: = \:  a^{\eps}_0(\zeta, \xi_n) \:
\prod_{j=1}^r \,
a^{\eps}_j\left(\zeta, \eps^{\gamma_j} \xi_n \right), 
\end{equation}
where 
$$ a^{\eps}_0(\zeta, \xi_n) \: := \: C^{\eps}_{\chi} \, 
\prod_{z \in {\cal Z}^{\eps}_0(\zeta)} 
\, \left( \xi_n \, - \, z  \right), $$
$$ a^{\eps}_{j}(\zeta, \eta) \: := \:\prod_{z \in {\cal Z}^{\eps}_j(\zeta)}  \, \left( 
 \eta \, - \, \eps^{\gamma_j} z  \right), $$
for $1 \le j \le r$. Note that for all $\zeta$ in ${\cal W}$,  
and for $0 < \eps \le \eps_0^{\zeta}$, we have 
$$ a^{\eps}_j(\zeta, \cdot) \: = \: a^{\eps, \zeta}_j $$
where $a^{\eps, \zeta}_j$ was introduced in \eqref{factor}. We can extend
 $a^{\eps}_j$ to $\eps = 0$ by 
$$ a^0_0(\zeta, \cdot) \: := \: \prod_{i=s+1}^m \left( \cdot \, - \,
\xi_{i,0}^{\zeta} \right), \quad   a^0_j(\zeta, \cdot) \: := \:
\prod_{p_i = \gamma_j} \left(  \cdot \, - \, \eta_{i,0}^{\zeta} \right), $$
for $1 \le j \le r$. We emphasize that all $a^0_j$ are  globally defined in
$\zeta$,
{\it i.e.} on $\{ | \zeta | > R \}$. On the contrary, the 
$a^{\eps}_j$ are only locally defined in $(\eps, \zeta)$, on an open
subset  $\displaystyle 
[0, \eps_0(\overline{\zeta})) \times {\cal
W}(\overline{\zeta})$. The regularity of these various
polynomials is given in 
\begin{prop} \label{regularity}
For all $0 \le j \le r$:
\begin{description}
\item[i)] The coefficients of  $a^0_j$ are symbols in $\zeta$, 
(restricted to  $|\zeta| > R$). 
\item[ii)] Let $|\overline{\zeta}| > R$. The
coefficients of  $a^{\eps}_j$ are smooth
functions of $(\eps, \zeta)$ in a neighborhood of 
$(\eps = 0, \overline{\zeta})$. 
\end{description}
\end{prop}
 {\bf Proof:} 
We use the notations introduced in the proof of proposition
\ref{zeros}. Again, we deal only with $j = r, r-1$, and leave the
general induction argument to the reader. 

\medskip
{\bf i)} In the course of the  previous proof, we have shown that the
coefficients of 
$$ a^0_{r, \chi}(\zeta, \eta) \: = \: C_{\chi} \, 
\eta^{m-m_r} \, \prod_{p_i = \gamma_r}
\, \left( \eta \, - \, 
\eta_{i,0}^{\zeta} \right) \: = \:
  C_{\chi} \, \eta^{m-m_r} \, a^0_r(\zeta, \eta)  $$ 
are polynomial in $\zeta$, $| \zeta | > R$, so that 
 the coefficients of $a^0_r$ share the same property.
Also, the coefficients of 
\begin{align*} 
a^0_{r-1,\chi}(\zeta, \eta) \: & = \: C_{\chi} \,
 \eta^{m-m_r-m_{r-1}}  \, \prod_{p_i =
\gamma_{r-1}} \left( \eta \, - \, \eta_{i,0}^{\zeta} \right) \, \left( 
 \prod_{p_i = \gamma_r}  -\eta_{i,0}^{\zeta} \right). 
 \\
&  \: = \: C_{\chi} \, \eta^{m-m_r-m_{r-1}} \, a^0_{r-1}(\zeta, \eta) 
 \left( 
 \prod_{p_i = \gamma_r}  -\eta_{i,0}^{\zeta} \right). 
\end{align*}
are polynomial in $\zeta$. As 
$$C_{\chi} \,  \left( \prod_{p_i = \gamma_r}  - \eta_{i,0}^{\zeta}
\right) $$
 is  the  coefficient of  order $m - m_r$ of $a^0_{r, \chi}$, it is
 a smooth symbol, which  does not cancel by
definition of the $\eta_{i,0}^{\zeta}$. We deduce that the coefficients
of  $a^0_{r-1}(\zeta, \cdot)$ are smooth symbols in $\zeta$ as well.

\medskip
{\bf ii)} In the course of previous proof, we have shown the existence of
a curve $\Gamma_r$ (independent on $\zeta \in {\cal W}$, $\eps \le
\eps_0 )$, enclosing the set $\displaystyle 
\eps^{\gamma_r} \, {\cal Z}^{\eps}_r(\zeta)$, and none of the other
roots of $a^{\eps}_{r, \chi}(\zeta, \cdot)$. By Cauchy's formula, we
deduce for all $k \in \N$,     
$$ \sum_{\eta \, \in \, \eps^{\gamma_{r-1}} \,  {\cal Z}^{\eps}_r} \, \eta^k 
\: = \: \frac{1}{2 i
\pi} \, \int_{\Gamma_r} \, \tilde{\eta}^k \, \frac{\pa_{\tilde{\eta}} \,
a^{\eps}_{r, \chi}(\cdot, \tilde{\eta})}{a^{\eps}_{r,\chi}(\cdot,
 \tilde{\eta})}\, d\tilde{\eta}. $$
The right-hand side, so the left-hand side, defines a smooth function of
$(\eps, \zeta)$ in $[0, \eps_0) \times {\cal W}$. Now, it is well known 
that  the coefficients of a polynomial are themselves  polynomial in  such  symmetric functions of the roots. 
 We deduce that the coefficients of $a^{\eps}_r(\zeta, \cdot)$ are
equally  smooth. 
We can proceed similarly for $a^{\eps}_{r-1}$, through the formula 
$$ \sum_{\eta \,  \in \,  \eps^{\gamma_{r-1}} \, {\cal Z}^{\eps}_{r-1}} \, \eta^k
 \: = \: \frac{1}{2 i
\pi} \, \int_{\Gamma_{r-1}} \, \tilde{\eta}^k \, \frac{\pa_{\tilde{\eta}} \,
a^{\eps}_{r-1, \chi}(\cdot, \tilde{\eta})}{a^{\eps}_{r-1, \chi}(\cdot,
\tilde{\eta})} \, d\tilde{\eta}, $$
and the smoothness of $a^{\eps}_{r-1}$ follows as well. $\Box$

\medskip
On the basis of (H1), we hope for an asymptotics of type
\eqref{asymptotique} with $\alpha_j = \eps^{\gamma_j}$ for some indices
$j$. Nevertheless, for such an asymptotics to be true, it is reasonable that
$r$, $(\gamma_1, m_1)$, \ldots, $(\gamma_r, m_r)$ depend neither on the
local chart $\chi$, nor on $y \in Y$.  With above notations, we have :
 for all $| \zeta | > R$, for all $1 \le j \le r$, 
\begin{equation} \label{a0j} 
a^0_j(\zeta, 0) \: = \:  (-1)^{m_j} \, \prod_{p_i = \gamma_j}
\eta_{i,0}^{\zeta} \neq 0. 
\end{equation}
By proposition \ref{regularity}, $a^0_j(\zeta, 0)$ is a symbol in 
$\zeta$ (restricted to $| \zeta | > R$). By equation
\ref{a0j}, it does not cancel. In general, this property of no
cancellation is not preserved by a change of variable. Therefore,
the assumption (H1), as well as $r$,$(\gamma_1, m_1)$, \ldots,
$(\gamma_r, m_r)$ depend on the local charts. To overcome this problem, we
need to make  a stronger assumption: 

\mspace
{\bf (H2)  For all $1 \le j \le r$, $a^0_j(\cdot, 0)$ is elliptic.}

\medskip
Note that, up to consider  a larger
  $R$, (H2) implies \eqref{a0j}. Note also that (H2) is
preserved by a change of the tangential variable $x'$.

\subsection{Asymptotic invariance and assumption (H3)}
 We can now state the following invariance result:
\begin{theo} {\bf (Invariance through diffeomorphism)} \label{invariance}

\sspace
Assume (H1)-(H2). Let $\left({\cal O}_{\Psi}, \Psi\right)$ a
local chart in $X$ around $y$ such that 
$$ \Psi_n(y) = 0, \quad \frac{\pa \Psi_n}{\pa t }(y) \neq 0. $$
Let (H1)$_{\psi}$-(H2)$_{\psi}$ be the same as 
(H1)-(H2), with $\Psi$ in place of $\chi$. Then,
(H1)$_{\psi}$-(H2)$_{\psi}\displaystyle$  holds, with the same   $r$,
$(\gamma_1, m_1)$, \ldots, $(\gamma_r, m_r)$.
\end{theo}
\begin{rem}
In short, theorem \ref{invariance} states that the conjunction 
(H1)-(H2) is intrinsic, as well as the main features of the
singularities (exponents $\gamma_j$ and multiplicities $r, m_j$). 
We stress that (H1)-(H2) involves the whole symbol 
$a^{\eps}_{\chi}(y, \cdot)$ and not only its principal
symbol. Therefore, it is not obvious that it should be preserved by another
choice of coordinates. 
\end{rem}

\medskip
{\bf Proof:} 
Let us first set a notation: 
for any smooth $\phi = (\phi', \phi_n) : U  \rightarrow
V$, $0 \in U$, for any symbol $P(\zeta, \eta)$ defined on $\R^{n-1} \times
\R$, and for any $\gamma \in \R$, we introduce
\begin{equation*}
\phi_*^{\gamma} P(\zeta, \eta) \, := \,
e^{i \phi(0) \cdot \left(\zeta, \eps^{-\gamma} \, \eta\right)} \left[ 
P(D_{x'}, \eps^\gamma D_{x_n}) \, 
 e^{-i \phi(x) \cdot \left(\zeta, \eps^{-\gamma} \, \eta\right)}
\right]_{\vert \, x = 0}
\end{equation*}
We also denote 
$$ \phi_* P(\zeta, \eta) \: := \: \phi_*^0 P(\zeta, \eta). $$

\medskip
Let $({\cal O}, \chi)$ and $({\cal O}_{\Psi}, \Psi)$ the local charts of the
theorem.  It is clear that (H1) and (H2) (and $r$,
$(\gamma_1, m_1)$, \ldots, $(\gamma_r, m_r)$) are not affected by a translation
 of each coordinate. Therefore, with no loss of generality,
we can assume that
$\displaystyle  \chi(y) = \Psi(y) = 0$. 
Let $\phi = \Psi \circ \chi^{-1}$. It is a smooth diffeomorphism between
two neighborhoods $U, V$ of $0$.
Let now $\displaystyle 
P^{\eps}(\zeta, \eta)$ a polynomial in $\eps, \zeta, \eta$. 
 The proof
of theorem \ref{invariance} will follow from the study of $\phi_*^{\gamma}
P^{\eps}$. Precisely, we will show:
\begin{prop} \label{diffeo}
Let $\gamma > 0$. Then, uniformly in each compact subset of 
$\R^{n-1} \times \R$, 
\begin{equation}  \label{limite}
\lim_{\eps \rightarrow 0} \phi_*^{\gamma}  \: P^{\eps}(\zeta, \eta) \:  = \:
\overline{\phi}_* \: P^0(\zeta, \eta), 
\end{equation}
where $\displaystyle \overline{\phi}(x',x_n) \: = \: \left( \phi'(x',0), \:
\pa_{x_n}\phi_n(x', 0) \, x_n \right)$. 
\end{prop}
\begin{rem}
If we denote
$$ \sum_{j} P^0_j(\zeta) \, \eta^j$$
 the expansion of $P^0$ with respect to $\eta$, then 
$$ \overline{\phi}_* \: P^0(\zeta, \eta) \: = \: \sum_{j} Q^0_j(\zeta) 
\, \eta^j,$$
with 
\begin{equation}
Q^0_j(\zeta)  \: = \: \left[ P^0_j(
D_{x'}) \left(  e^{-i \phi'(x', 0) \cdot \zeta} \: \frac{\pa
\phi_n}{\pa x_n}(x',0)^j \right) \right]_{\vert \,
x' = 0}. 
\end{equation}
The symbols $Q^0_j$ are polynomial in $\zeta$. More precisely, classical
symbolic calculus (see \cite{AlGe})
 provides the following (finite) expansion:  
\begin{equation} \label{classicexp}
 Q^0_j(\zeta) \: \sim \: 
\sum_{\alpha \ge 0} \frac{1}{\alpha!}  \, \pa^\alpha_{x'} \left[ 
e^{i r'\left(x'\right)} \: \frac{\pa
\phi_n}{\pa x_n}(x', 0)^j \right]_{\vert x' = 0} \:
D^\alpha_\zeta P^0_j\left(d_{x'}\phi'(0) \zeta \right), 
\end{equation}
 where
$$r'\left(x'\right) \: := \: \phi'\left(x', 0\right) 
 \: - \: d_{x'}\phi'(0) \, \left( x'\right). $$
The terms of rank $\alpha$ in this expansion are polynomial in $\zeta$, 
 of degree less than $\mbox{deg}\left(P^0_j\right) - |\alpha|/2$. (see
again \cite{AlGe} for details).
\end{rem}
\begin{rem} 
In the case $\gamma = 0$, the limit \eqref{limite} is easily replaced by 
\begin{equation} \label{limite0}
\lim_{\eps \rightarrow 0} \phi_*  P^{\eps}(\zeta, \eta)  \: =
\: \phi_* P^0(\zeta, \eta).
\end{equation}
\end{rem} 
We postpone temporarily the proof of the proposition, and we show how to
deduce theorem \ref{invariance} from it. We must prove that 
(H1)$_{\Psi}$-(H2)$_{\Psi}$ holds, that is we can replace the local symbol
$a^\eps_\chi$ by the other local symbol $a^\eps_\Psi$. With a slight abuse
of notation, we identify $x$ and $\chi(x)$ so that we work with $x \in U$,
$U$ a neighborhood of $y=0$ in $\R^n$. Hence,
\begin{equation} \label{changeoperateur}
 a^\eps_\Psi(0, \cdot) \: = \: \phi_* \,  a^\eps_\chi(0, \cdot). 
\end{equation}

\medskip
{\em Proof of} (H1)$_{\Psi}$: We consider 
$$ a^{\eps}_{r, \Psi}(\zeta, \eta) : = \eps^{-\beta_r} 
a^{\eps}_{\Psi}(0, \zeta, \frac{\eta}{\eps^{\gamma_r}}),$$
and, as in previous proofs,
$$ a^{\eps}_{r, \chi}(\zeta, \eta) : = \eps^{-\beta_r} \, 
a^{\eps}_{\chi}(0, \zeta, \frac{\eta}{\eps^{\gamma_r}}).$$
 Equation \eqref{changeoperateur} becomes
$$  a^{\eps}_{r, \Psi} \: = \: \phi_*^{\gamma_r} \, a^{\eps}_{r, \chi}. $$
We have seen in the previous proofs that   $a^{\eps}_{r, \chi}$ is 
polynomial in $\eps^{\alpha}, \zeta, \eta$ for some appropriate 
$\alpha > 0$. Moreover,  
\begin{equation} \label{expressiona0r}
a^{0}_{r, \chi}(\zeta, \eta) \: =\: C_{\chi} \,
\eta^{m-m_r} a^0_r(\zeta, \eta), \quad |\zeta| > R. 
\end{equation}
 If we expand 
$$ a^0_{r, \chi}(\zeta, \eta) \: = \:  \sum_{j}   
a_{r,j}^0(\zeta) \: \eta^j, $$ 
we get from proposition \ref{diffeo} that, for $|\zeta|$ large enough, 
$$ \lim_{\eps \rightarrow 0} \:   a^{\eps}_{r, \Psi}(\zeta, \eta) \: = \: 
\sum_j  b_{r,j}^0(\zeta) \: \eta^{j}, $$
with 
$$  b_{r,j}^0(\zeta) =    
\left[ a^0_{r,j}(D_{x'}) 
\left(  e^{-i \phi'_0(x') \cdot \zeta} \: \frac{\pa
\phi_n}{\pa x_n}(x',0)^j \right) \right]_{\vert \,
x' = 0} $$ 
Using \eqref{expressiona0r},  we deduce that 
the term of lowest degree (with respect to
$\eta$) in  $\displaystyle \lim_{\eps \rightarrow 0} \: 
a^{\eps}_{r, \Psi}(\zeta, \eta) \:$ is $\: \displaystyle 
b_{r, m-m_r}^0(\zeta) \:
\eta^{m-m_r}$. Moreover, by \eqref{classicexp}, the principal symbol of
$\displaystyle b_{r,m-m_r}^0$ is    
\begin{align*}
 \sigma\left(b_{r, m-m_r}^0\right)(\zeta, \eta) \: & =  \:
\sigma\left(a^0_{r,m-m_r}\right)(d\phi'_0(0) \, \zeta) \\
& = C_{\chi} \:
 \sigma\left(a^0_r(\cdot, 0)\right)(d\phi'_0(0) \, \zeta). 
\end{align*}
Note that $d_{x'}\phi'(0)$ is a  diffeomorphism  in
$\R^{n-1}$. By assumption (H2), $a^0_r(\cdot, 0)$ is elliptic, hence
$\displaystyle b_{r, m-m_r}^0$ is still elliptic. As a consequence, it
does not vanish  for $|\zeta|$ large enough. Therefore, at fixed $\zeta$, the
polynomial $\lim_{\eps \rightarrow 0} \: a^{\eps}_{r, \Psi}(\zeta, \cdot)$
has $m_r$ roots away from zero. By Rouch\'e's theorem, it is still the case
for $\eps \le \eps_0$ small enough. Back to the original variables, we find
$m_r$ roots of  $a^{\eps}_{\Psi}(0, \zeta, \cdot)$ satisfying 
$$ \frac{\delta}{\eps^{\gamma_r}} \: \le \: | \xi_n | \le
\frac{\delta^{-1}}{\eps^{\gamma_r}}. $$

\medskip
With similar arguments, we can find $m_{r-1}$ roots of $a^{\eps}_{\Psi}(0,
\zeta, \cdot)$ which scale like $\eps^{-\gamma_{r-1}}$. We consider this time
$$ a^{\eps}_{r-1,\Psi} \: :=  \: \eps^{-\beta_{r-1}} \, 
a^{\eps}_{\Psi}(0, \zeta,
\frac{\eta}{\eps^{\gamma_{r-1}}}), $$
as well as $a^{\eps}_{r-1,\chi}$. Again, $a^{\eps}_{r-1,\chi}$ is a
polynomial in $(\eps^{\alpha}, \zeta, \eta)$ and 
$$ a^{0}_{r-1, \chi}(\zeta, \eta) \: =\: C_{\chi} \,
\eta^{m-m_r-m_{r-1}} \: a^0_{r-1}(\zeta, \eta) \: a^0_r(\zeta, \eta).$$
We then ap\-ply the same rea\-so\-ning as above, replacing  $\dis a^0_r(\cdot, 0)$ 
by the product $\dis a^0_{r-1}(\cdot, 0) \: a^0_r(\cdot, 0)$, 
which is still elliptic.

\medskip
Proceeding recursively, we find, for all $1 \le j \le r$, $m_j$ roots of
$a^{\eps}_{\Psi}(0, \zeta, \cdot)$ such that 
$$ \frac{\delta}{\eps^{\gamma_j}} \: \le \: | \xi_n | \le
\frac{\delta^{-1}}{\eps^{\gamma_j}},$$
for $\eps \le \eps_0$ small enough, $\delta = \delta(\eps_0)$ small
enough. Finally, by \eqref{limite0}, 
$$ \lim_{\eps_{\rightarrow 0}} \:  a^{\eps}_{\Psi}(0, \zeta, \xi_n) \: = \: 
\phi_* \left( a^0_0(\zeta, \xi_n) \, \prod_{j=1}^r a^0_j(\zeta, 0)
\right). $$
At fixed $\zeta$, the right-hand side defines a polynomial in $\xi_n$,
of degree $m-s$. By Rouch\'e's theorem, this yields $m-s$ roots of 
$a^{\eps}_{\Psi}(0, \zeta, \cdot)$ with 
$$ 0 \: \le \: | \xi_n | \: \le \delta^{-1}. $$ 
As the degree of $a^{\eps}_{\Psi}(0, \zeta, \cdot)$ is $m = m-s + \sum
m_j$, we obtain in this way all the roots. Assumption (H1)$_{\Psi}$
follows with the same $r$, $\gamma_1$,\ldots,  $\gamma_r$. 

\medskip
{\em Proof of} (H2)$_{\Psi}$: We have just established (H1)$_{\Psi}$. Just
as for $a^{\eps}_\chi$, we can factorize 
\begin{equation} \label{factor3}
a^{\eps}_{\Psi}(y, \zeta, \xi_n) \: = \: \tilde{a}^\eps_0(\zeta, \xi_n)
\: \prod_{j=1}^r \tilde{a}^\eps_j(\zeta, \eps^{\gamma_j} \xi_n) 
\end{equation}
(see factorization \eqref{factor2}), and we have to show that: for all $1
\le j \le r$, $\tilde{a}^0_j(\zeta, 0)$ is elliptic in a neighborhood of
$y$. But  from above identities, we have got easily:
$$ C_{\Psi} \: \sigma\left(\tilde{a}^0_r(\cdot, 0)\right)(\zeta) \: = \: 
 C_{\chi} \: \sigma\left(a^0_r(\cdot, 0)\right)(d\phi'_0(0) \zeta), $$
and
\begin{multline*}
 C_{\Psi} \:\: \sigma\left(\tilde{a}^0_{r-1}(\cdot, 0)\right)(\zeta) \:\: 
 \sigma\left(\tilde{a}^0_r(\cdot, 0)\right)(\zeta)  \: =  \\  
 \: C_{\chi} \:\: \sigma\left(a^0_{r-1}(\cdot, 0)\right)( d_{x'}\phi'(0) \,
 \zeta)  \:\: \sigma\left(a^0_r(\cdot, 0)\right)(d_{x'}\phi'(0)\, \zeta), 
\end{multline*}
and so on. As ellipticity is preserved by inversion and product, we deduce 
(H2)$_{\Psi}$ from a simple recursion. $\Box$

\medskip
{\bf Proof of proposition \ref{diffeo}}: The proof relies on the
representation of $\phi_*^\gamma P^\eps$ as an oscillatory integral. The
asymptotic analysis of such integral will be performed through ``stationary
phase type'' theorems. This approach is very classical in symbolic calculus, to
establish the properties of conjugation, product, or transformation under a
diffeomorphism (see, among many, textbooks \cite{ChPi, AlGe}). 
However, in the standard
setting, one is mainly concerned with the principal symbol, so that the
small parameter in the asymptotics is the wavelength $1/|\xi|$. In our
framework, the natural idea is to work with the parameter $\eps$ instead of
the wavelength. 

\medskip
Let $u \in {\cal C}^{\infty}_c(U)$, $\zeta \in
\R^{n-1}$, $\eta \neq 0$. We wish to compute
\begin{equation*} 
I^\eps  \: = \: 
\left[ P^\eps(D_{x'}, \eps^\gamma D_{x_n}) \left( e^{-i
\phi(x) \cdot (\zeta, \eps^{-\gamma} \eta)} u(x) \right)
\right]_{\vert \, x = 0}  
\end{equation*}
$I^\eps = I^\eps(\zeta, \eta)$  is defined by the oscillatory
integral:
\begin{equation*}
I^{\eps} \, := \, \frac{1}{(2 \pi)^n} \int e^{-i x \cdot
(\tilde{\zeta}, \tilde{\eta})} \: P^{\eps}\left(\tilde{\zeta},
\eps^\gamma \tilde{\eta}\right) \: u(x) \: 
 e^{i \phi(x) \cdot (\zeta,  
\eps^{-\gamma} \eta)}  
\: dx \, d(\tilde{\zeta}, \tilde{\eta}). 
\end{equation*}
Following the scheme of
\cite[page ??]{AlGe}, we make the change of variables:
$$\tilde{\eta} \: := 
\: \tilde{\eta} \: -  \: \frac{\pa_{x_n} \phi_n(x',0)}{\eps^\gamma} 
\eta. $$
By the Fubini theorem for oscillatory integrals, we can write
$$ I^\eps \: = \: \frac{1}{(2 \pi)^n} \int \, e^{-i x' \cdot
\tilde{\zeta}} \: J^\eps(x',\tilde{\zeta} ) \: dx' \,
d\tilde{\zeta}, $$ 
with
\begin{multline*}
J^\eps(x',\tilde{\zeta}) := \int e^{-i x_n \, \tilde{\eta}}
 \: e^{i \, \eps^{-\gamma} \,  r_{\eta}(x)}
 \: u(x) \:
e^{i \phi'(x) \cdot \zeta} \\
P^\eps\left(\tilde{\zeta}, \eps^\gamma \tilde{\eta} + \pa_{x_n}
\phi_n(x', 0) \eta \right) \: dx_n d\tilde{\eta},
\end{multline*}
where  
$$ r_{\eta}(x) \: := \:  \eta \, \left( 
\phi_n(x)  - \pa_{x_n} \phi_n(x', 0)\right).$$
 We set 
$$ \Phi^{\eps}(x', \tilde{\eta}) := \int e^{-i x_n \,
\tilde{\eta}} \:  e^{i \frac{r_\eta(x)}{\eps^\gamma}} 
\: u(x) \: e^{i \phi'(x) \cdot \zeta} \:
dx_n, $$
so that 
$$ J^\eps(x',\tilde{\zeta}) \: = \: \int 
P^\eps\left(\tilde{\zeta}, \, \eps^\gamma \tilde{\eta} \, + \,
 \pa_{x_n} \phi_n(x', 0) \, \eta \right) 
\:  \Phi^{\eps}(x',
\tilde{\eta}) \: d\tilde{\eta}. $$
For fixed $\eta, x'$, we introduce the phase
$$ \eps^{-\gamma} \, f(x_n) \: := \:\eps^{-\gamma} \left( 
  r_{\eta}(x) \: - 
\: \eps^\gamma \tilde{\eta}  \, x_n \right), $$
with derivative
$$ f'(x_n) \: = \: \eta \left( \pa_{x_n} \phi_n(x) -  
\pa_{x_n}  \phi_n(x', 0)\right) - 
\eps^\gamma \tilde{\eta}. $$

\medskip
\begin{itemize}
\item  If $\: \left| \, \eta \, \pa_{x_n} 
\phi_n(x', 0) \, + \, \eps^\gamma
\tilde{\eta} \, \right| \: < \:  1/C $, $C$ such that $| \, \eta
\,  \pa_{x_n} \phi(x) \, 
| \: \ge \:  2/C$ for $x$ in a neighborhood of $0$, then $|
f'(x_n)| \: \ge \: 1/C$. 
\item If $\: \left|  \, \eta \, \pa_{x_n} 
\phi_n(x', 0) \, + \, \eps^\gamma
\tilde{\eta}\,  \right| \: \ge \: C$, $C$ such that $|  \,\eta 
\pa_{x_n} \phi(x) \,
| \:\le\: C/2$ for $x$ in a neighborhood of $0$, then $|
\,f'(x_n)\,| \: \ge \: C/2$.
\end{itemize}
We deduce from the non-stationary phase theorem (see \cite[page ??]{AlGe}) 
that for large enough $C$, outside
$$ \frac{1}{C} \: \le \: \left| \eta  \pa_{x_n} 
\phi_n(x', 0) + \eps^\gamma
\tilde{\eta} \right| \: \le \: C, $$
 for all $k$,
\begin{equation} \label{decroissance}
 \left|  \Phi^{\eps}(x',
\tilde{\eta}) \right| \: \le \: C_k \, \left( 1 + |\tilde{\eta}| +
\frac{\eta}{\eps^\gamma} \right)^{-k}.   
\end{equation}

\medskip
Still following \cite{AlGe}, we introduce a truncation function $\alpha$ such 
that $\alpha \in {\cal C}^{\infty}_c(\R)$, $\alpha(\xi_n) = 1$ for $C^{-1}
\le |\xi_n| \le C$, $\alpha(\xi_n) = 0$ for $|\xi_n| \le (2C)^{-1}$. We
divide
$$ J^\eps(x', \tilde{\zeta}) = J^\eps_1(x', \tilde{\zeta}) + 
J^\eps_2(x', \tilde{\zeta}) $$
with 
\begin{equation*}
 J^\eps_2 \: = \: \int \alpha\left( \pa_{x_n}\phi_n(x', 0) +
 \frac{\eps^\gamma}{\eta} \tilde{\eta} \right)   
P^\eps\left(x, \, \tilde{\zeta}, \, \eps^\gamma \tilde{\eta} \, + \, 
 \pa_{x_n} \phi_n(x',x_n) \, \eta \right) 
\: \Phi^{\eps}(x',\tilde{\eta}) \: d\tilde{\eta}. 
\end{equation*}
Note that $J^\eps_1$ and $J^\eps_2$ are polynomials in $\tilde{\zeta}$ 
 (as the only dependence on $\tilde{\zeta}$ is in $P^\eps$). From the
 estimate \eqref{decroissance}, we deduce
$$ \left| \pa_{(x', \tilde{\zeta})}^\beta \,  J^\eps_1(x',
\tilde{\zeta}) \right| \: \le \: C_{\beta,k} \, 
|\tilde{\zeta}|^{N - \beta} \, \eps^k, \quad N \in \N, \quad \forall \beta,
 k.$$
Hence, 
$$I^\eps_1 \: := \:  \frac{1}{(2 \pi)^n} \int e^{-i x' \cdot \tilde{\zeta}} 
\; J^\eps_1( x',\tilde{\zeta}) \: dx' \,  d\tilde{\zeta} $$
satisfies $I^\eps_1 = O(\eps^k)$ for all $k$. Moreover, a rapid look at the
proof shows that all estimates hold uniformly in every compact subset of
$\R^{n-1} \times \R_*$. 

\medskip
It remains to estimate the ``stationary terms'' $J^\eps_2$ and $I^\eps_2$. 
For this, we will use a lemma whose statement and 
proof can be found in \cite{AlGe}:
\begin{lem} (\cite{AlGe})  \label{alinhac}

\sspace
Let $r= r(z,\theta) \in {\cal C}^{\infty}(\R \times {\cal O} )$, 
${\cal O}$  open subset of $\R^p$, 
 such that $\pa_z r(0, \theta) = 0$.  Let $b^\lambda = 
b^\lambda(z, \eta, \theta) \in {\cal C}^{\infty}(R \times \R \times 
{\cal O} )$,
depending on a parameter $\lambda > 0$ so that, for some $M$,
\begin{equation} \label{decreaseb}
 \forall \, \beta, \beta', \quad \left| \pa_z^\beta 
\pa_\eta^{\beta'}
b^\lambda(z,\eta, \theta, \lambda) \right| \le C_{\beta,\beta'} \:
\lambda^{M - \beta'}.
\end{equation}
Assume moreover that $b^{\lambda}$ has compact support in $z$, uniformly
with respect to $\eta$, $\theta$, $\lambda$. Then, the integral
$$ J(\lambda) \: = \:  (2 \pi)^{-n} \int e^{-i \, z \, \eta} e^{i \lambda r(z,
\theta)} \: b^\lambda(z, \eta, \theta) \: dz \, d\eta $$
has the following asymptotic expansion for $\lambda \rightarrow +\infty$
$$  J(\lambda) \: \sim \: \sum_{\alpha \ge 0} \frac{1}{\alpha!}
\pa_z^\alpha \left( e^{i \lambda r(z, \theta)}D^\alpha_\eta
b^{\lambda}(z,0,\theta) \right)\vert_{z=0}, $$
where terms of rank $\alpha$ are bounded by $\lambda^{M-\alpha/2}$.
This expansion holds locally uniformly with respect to $\theta$. 
\end{lem}
The function $J_2^\eps(x, \tilde{\zeta})$ 
is polynomial in $\tilde{\zeta}$, with
coefficients that read
$$
\int  e^{-i x_n \, \tilde{\eta}} 
 \: e^{i \eps^{-\gamma} r_{\eta}(x) \, \eta} \:
 b_j^\eps(x_n,
 \tilde{\eta}, \zeta, \eta, x') \: dx_n.$$
Such integrals can be easily put in the framework of lemma \ref{alinhac}, with
 $\displaystyle 
\lambda :=  \eps^{-\gamma}$, $z:=x_n$, $\eta := \displaystyle 
\tilde{\eta}$, and 
 $\displaystyle \theta := (\zeta, \eta, x')$. Thus, 
$r(z, \theta) := r_{\eta}(x)$ satisfies $\pa_z r(0, \theta) =
 0$. Besides, the function 
$$b^{\lambda}(z,\eta,\theta) :=
 b_j^\eps(x_n, \tilde{\eta}, \zeta, \eta, x')$$ 
satisfies  estimate \eqref{decreaseb} with $M =0$. Its support in $z$ is
contained in the support of $u$, 
so uniform with respect to the other variables. Hence, we can apply lemma
\ref{alinhac}. The first term of the expansion provides
\begin{equation*}
 J_2^\eps(x, \tilde{\zeta}) \xrightarrow[\eps \rightarrow 0]{} 
J_2^0(x, \tilde{\zeta}) \: := \: 
e^{i \, \phi'(x', 0)} \:
P^0\left(\tilde{\zeta}, \, \pa_{x_n} \phi_n(x', 0) \,
\eta\right)  \: u(x', 0)   
\end{equation*}
and consequently
\begin{equation*}
I_2^\eps \: \xrightarrow[\eps \rightarrow 0]{} \: \frac{1}{(2 \pi)^n}\int 
e^{-i x' \cdot \tilde{\zeta}} \, J_2^0(x, \tilde{\zeta}) \: 
dx' \, d\tilde{\zeta}
\end{equation*}
with convergence in every compact subset of $\R^{n-1} \times
\R_*$. As we handle polynomials in $\eta$, the uniform convergence extends
to compact subsets of $\R^{n-1} \times
\R$. By taking $u$  equal to $1$ in a neighborhood of $0$ in $U$, we obtain
the result.  $\Box$





\medskip
We end this section on asymptotic invariance with a definition: 
\begin{defi}
Let $y \in Y$. We say that $a^\eps$ is uniformly singular at $y$ if
(H1)-(H2) holds at $y$. The corresponding $\gamma_j$'s, $1 \le j \le r$ are
called singular exponents (of $a^\eps$ ) at $y$, and the $m_j$'s are their
multiplicity.
\end{defi} 
These definitions make sense because of theorem
\ref{invariance}. We recall that the singular exponents are positive
rational numbers (see corollary \ref{holom2}). 
\begin{defi}
Let $Y' \subset Y$. We say that $a^\eps$ is uniformly singular on $Y'$, if
it is uniformly singular at all $y \subset Y'$, with constant singular
exponents and multiplicities. 
\end{defi}
In the whole sequel, we will make the following assumption:

\mspace
{\bf (H3)  The operator  $a^\eps$ is uniformly singular on each connected
component of the boundary  $Y$.}




\section{Boundary layer sizes  and equations} \label{blsizes}
\subsection{Definitions and hypothesis (H4)}
Let $Y'$ a connected component of $Y$. By assumption (H3), $a^\eps$ is
uniformly singular on $Y'$. We  denote by  $\gamma_1$ to
$\gamma_r$ the singular exponents, and $m_1$ to  $m_r$ their
respective multiplicity.
In a region $Y' \times (0,T')$, equation \eqref{elliptique} reads 
$$ a^\eps(x, D_y, D_t) u^\eps \: = \: 0. $$
Let $y \in Y'$. Let $({\cal O}', \chi')$ a local chart in $Y'$ around
$y$. The application 
$$ \chi : {\cal O}'\times (-T', T') \mapsto \R^n, \quad (y', t) \mapsto
(\chi'(y'), t) $$
defines a local chart in $X$ around $y$. Hence, we can define a local
operator $a^\eps_\chi(x, D_{x'}, D_t)$, with symbol $a^\eps_\chi(x, \zeta,
\xi_n)$. 

\medskip
On one hand, we have a finite expansion of the type 
\begin{equation} \label{akk}
 a^\eps\left( y, D_y, \frac{D_\theta}{\eps^{\gamma_j}}\right) \: = \:
\sum_{k=-k_0}^{k_0} \: \sum_{k'} \: 
A_{\substack{k,k'\\ \chi'}}(y, D_y) \, D_\theta^{k'}
\eps^{\alpha k} 
\end{equation}
for some appropriate $\alpha > 0$ and smooth differential operators
$A_{k,k'}$  on $Y'$.
 Hence, 
\begin{equation} \label{localsymbol1}
a^\eps_\chi\left(y, \zeta, \frac{\eta}{\eps^{\gamma_j}}\right) \: = \: 
\sum_{k=-k_0}^{k_0} \: \sum_{k'} \: A_{k,k'}(y, \zeta) \, \eta^{k'}
\eps^{\alpha k}.
\end{equation}

\medskip
On the other hand, with notations introduced in \eqref{factor2}, we have
for  $|\zeta| > R(y)$, 
\begin{equation} \label{localsymbol2}
\lim_{\eps \rightarrow 0} \, \eps^{\beta_j} a^\eps_\chi\left(y, \zeta,
\frac{\eta}{\eps^{\gamma_j}}\right) \: = \: C_\chi \, \eta^{m - \sum_{j'
    \ge j} m_{j'}} \: a^0_j(\zeta, \eta) \: \prod_{j' > j} 
a^0_{j'}(\zeta, 0)
\end{equation}
From \eqref{localsymbol1} and \eqref{localsymbol2} we deduce
\begin{equation} \label{localsymbol3}
 \lim_{\eps \rightarrow 0} \, \eps^{\beta_j} a^\eps\left( y, D_y,
\frac{D_\theta}{\eps^{\gamma_j}}\right) \: = \: D_\theta^{m- \sum_{j' \ge
    j} m_j'} a_{\gamma_j}(y, D_y, D_\theta), 
\end{equation}
where $a_{\gamma_j}$ is a smooth differential operator. It satisfies, 
 for $|\zeta| > R(y)$, 
$$   a_{\gamma_j, \chi}(y, \zeta, \eta) \: = \: 
 a^0_j(\zeta, \eta) \: \prod_{j' > j} a^0_{j'}(\zeta, 0)
$$
For $1 \le j \le r$, we say that the  operator $a_{\gamma_j}$
is a singular operator (associated to the singular exponent $\gamma_j$). 

\medskip
Special attention will be paid to localized solutions of equation 
\begin{equation} \label{equationbl}
a_\gamma \, u \: = \: f, 
\end{equation}
where $a_\gamma$ is a singular operator, associated to the singular
exponent $\gamma$. Indeed, the functions $v^{i,k}_j$ in expansion
\eqref{ansatz3} will satisfy equations of type \eqref{equationbl}. 
Therefore, we introduce the Fr\'echet space
$$ {\cal X}  \: := \: {\cal S}\left( \R^+_t ; \, {\cal C}^\infty(Y')\right). $$
We also denote 
$$ {\cal X}_0  \: := \:  \left\{ u \in {\cal X}, \quad a_\gamma u =  0
\right\}. $$ 
We set the following definition 
\begin{defi}
We say that $\gamma$ is a boundary layer exponent (and $\eps^\gamma$ a
boundary layer size) on $Y'$, if  there exist subspaces ${\cal X}'_0 \neq
\{ 0 \}$, ${\cal X}'$ with:
\begin{description}
\item[i)] ${\cal X}'_0 \subset {\cal X}_0, \quad {\cal X}'_0 \subset
 {\cal X}' \subset {\cal X}$. 
\item[ii)]$ D_\theta {\cal X}' = {\cal X}', \quad D_y {\cal X}' 
\subset {\cal X}', \quad \theta {\cal X}' 
\subset {\cal X}', \quad f \, {\cal X}' 
\subset {\cal X}'$ for all smooth functions $f \in {\cal C}^\infty(Y')$.
\item[iii)] For all $f \in {\cal X}'$, \eqref{equationbl} has a solution
  in ${\cal X}'$.
\end{description}
\end{defi}

\begin{rem}
Note that identity $D_\theta {\cal X}' = {\cal X}'$ implies that ${\cal
X}'$ is preserved by both $D_\theta$ and its inverse $D_\theta^{-1}$ (which
is obviously defined on ${\cal X}$). Note also that ${\cal X}'_0$ and 
${\cal X}'$ are not uniquely determined. However, the subspaces $+ \,{\cal
X}'_0$ and $+ \, {\cal X}'$ where the algebraic sum is taken over all
subspaces   ${\cal X}'_0$ and  ${\cal X}'$ satisfying i)-iii), are 
well-defined, and are the largest subspaces satisfying i)-iii)
\end{rem}

\begin{rem}
The definition of a boundary layer exponent is at first sight dependent on
the diffeomorphism 
$$ \chi : {\cal T}' \approx Y' \times (-T',T'),  \quad x \approx (y, t), $$
where ${\cal T}'$ is a tubular neighborhood of $Y'$.  However, if 
$$ \psi : \tilde{{\cal T}}' \approx Y \times (-\tilde{T}',\tilde{T}'),  
\quad x \approx (\tilde{y}, \tilde{t}) $$ 
is another such diffeomorphism, we can link the singular operators
$a_{\gamma, \chi}$ and $a_{\gamma, \psi}$  following computations of
theorem \ref{invariance}. Indeed, we have 
$$ a_{\gamma, \psi}(y, D_y, D_t) \: = \: \overline{\phi}_* \, 
a_{\gamma, \chi}(y, D_{\tilde{y}}, D_{\tilde{t}}), $$
where $\phi = (\phi', \phi_n) = \psi \circ \chi^{-1} \in Y' \times (-T_{min},
T_{min})$, where $T_{min} \, :=  \, \min(T', \tilde{T}')$. , and  
$$  \overline{\phi} : (y', t') \: \mapsto \:  \left(\phi'(y', 0), \, 
\frac{\pa  \phi_n}{\pa t}(y', 0) \, t' \right). $$
Thus, solutions of \eqref{equationbl} with $a_\gamma = a_{\gamma, \chi}$
and $a_\gamma = a_{\gamma, \psi}$ are deduced from each other by the change
of variable $\overline{\phi}$ . As soon as 
$$ \frac{\pa  \phi_n}{\pa t} \: \ge \: C_0 > 0,$$
this change of variables preserves the properties of a boundary layer
exponent. 
\end{rem}

We end this section with  an assumption that will ensure the existence
 of a singular (boundary layer part) in \eqref{ansatz1}:

\mspace
{\bf (H4) Among the singular exponents given by (H3), there exists at least
one boundary layer exponent.}

\subsection{Solvability}
The verification of (H4) relies on solutions of 
\eqref{equationbl} in ${\cal X}$. In some particular cases, one can give
some effective criteria that ensure the solvability of this equation.

\subsubsection{Order zero operators}
For many physical systems, the coefficients of the singular operators are
just smooth functions 
$$ a_\gamma(y, D_y, D_\theta) \: = \:  
 a_\gamma(y, D_\theta).$$
Hence, \eqref{equationbl} resumes to a collection of linear ODE's with
constant coefficients, indexed by $y$ in $Y'$. Thus, to determine if
$\gamma$ is a boundary layer exponent is relatively easy. For instance,
\begin{prop} \label{solv}
Assume that there exists a closed curve 
$$\Gamma \subset \left\{ \Im(\eta) > 0 \right\}$$
 enclosing the set 
$$ \left\{ \eta, \quad \exists y, a_\gamma(y, \eta) = 0 \right\} \: \cap \:
\left\{ \Im(\eta)> 0 \right\}$$
of all roots with positive imaginary part. If this set is not empty, 
 then  $\gamma$ is a boundary layer exponent, otherwise it is not. 
\end{prop}
{\bf Proof: }By Rouch\'e's theorem and smoothness of $a_\gamma$, the
number $m^+(y) $ of roots of $a_\gamma(y, \cdot)$ lying inside
$\Gamma$ is locally constant, As $Y'$ is connected, it is constant.
 It is then well-known (see \cite{Agr}) that the
homogeneous space ${\cal X}_0$ is spanned by functions of the type
$$ f(y)\, \omega_j(y, \theta), \: f \in {\cal C}^\infty(Y'), \: 
\omega_j(y, \theta) \: = \: \int_{\Gamma} e^{i \,\eta \, \theta} \:
a^{-1}_\gamma(y, \eta) \:  \eta^j \: d\eta, \quad j \in \N. $$
Moreover, the dimension of 
$\displaystyle \mbox{ Vect }\left\{ \omega_j(y, \theta), \: j\in \N
\right\}$ is $m^+$.  The result follows, using ${\cal X}'_0 = {\cal X}_0$,
and  
$$ {\cal X}' \: = \: \mbox{ Vect } \left\{ f(y) \, 
\: \int_{\Gamma} e^{i \,\eta \, \theta} \:
a^{-l}_\gamma(y, \eta) \:  \eta^{j} \: d\eta \, , 
\quad f \in {\cal C}^\infty(Y'), \quad l \in \N, j
\in \Z \right\} $$  
$\Box$

\medskip
As shown in proposition \ref{solv}, the boundary layer is connected to
roots with positive imaginary part. This echoes the formal method used in
\cite{Ger4}, where roots $\xi_n^\eps$ with $\Im \, \xi_n^\eps \rightarrow
+\infty$ were considered.   
\subsubsection{Operators with constant coefficients}
Suppose that $Y'= \R^{n-1}$ (or $\T^{n-1}$) and $a_\gamma$ has constant
coefficients. One can use the Fourier transform over $\R^{n-1}$ (or
$T^{n-1}$), which turns \eqref{equationbl} into 
$$ a_\gamma(\zeta, D_\theta) \hat{u} \: = \: \hat{f}. $$ 
This is still a collection of linear ODE's, with constant coefficients,
indexed by $\zeta \in  \R^{n-1}$ (or $\T^{n-1}$). 

\subsubsection{General case and link with spectral problems}
In general the singular operator $a_\gamma(y, D_y, D_\theta)$ will be a
"true" differential operator. But in many cases it will
give rise to a (spectral) decomposition which will allow us to solve the problem. In this paragraph we will describe 
briefly the methodology without focusing on technical details.

Let us remind that a  solution of a problem of type
\begin{equation}\label{thierry1}
a_\gamma(y, D_y, D_\theta) \, u(y,\theta) \ =\ f(y,\theta)
\end{equation}
 that decreases as $\theta\to +\infty$ can be handled by Fourier-Laplace transform in $\theta$. Denoting $\tilde
 u(y,\xi)\ := \ \int_0^\infty 
 e^{-\xi\theta}u(y,\theta)d\theta$, (\ref{thierry1}) is equivalent to an equation
 of the form : 
 \begin{equation}\label{thierry2}
 a_\gamma(y, D_y, \xi) \, 
\tilde u(y,\xi)\ =\ \tilde f(y,\xi)+F(y,\xi) :=G(y,\xi)
 \end{equation}
 where $F$ involves boundary terms issued from  integrations by parts in 
the Laplace transform.
 Let us suppose that $A_\xi:=a_\gamma(y, D_y, \xi)$ admits a ( sometimes
 spectral) 
 decomposition weakly on $L^2$ of the form: 
 \begin{equation}\label{decompo}
 A_\xi=\int g(\lambda,\xi)|\varphi_{\xi,\lambda}><\varphi_{\xi,\lambda}|d\lambda
 \end{equation}
 where $|\varphi_{\xi,\lambda}><\varphi_{\xi,\lambda}|$ denotes the orthogonal projector on $\varphi_{\xi,\lambda}$, and the
 (spectral) parameter $\lambda$ is either complex or real. Let us 
 suppose moreover that the function $g(\lambda, \xi)$ and the family \
of vectors $\varphi_{\xi,\lambda}$ are analytic in $\xi$ (let us
 mention here that a naive definition of analytic vectors is enough for the 
examples we have in mind:  for a precise definition see
 \cite{Reed:1980}). We can solve (\ref{thierry2})  by using the decomposition 
(\ref{decompo}) . Indeed denoting 
 
 $$\hat G(\lambda,\xi) \, :=  \, 
<\varphi_{\xi,\lambda}|G(\lambda,.)>:=\int\overline{\varphi_{\xi,\lambda}}
G(\lambda,\xi)d\xi$$ 

we get:
\begin{equation}\label{thierry3}
u(y,\theta) = 
\frac 1 {2\pi i}\int d\lambda\int_{x-i\infty}^{x+i\infty}d\xi
\frac{1}{g(\lambda,\xi)}\hat G(\lambda,\xi)e^{\xi\theta}
\varphi_{\xi,\lambda}(y).
\end{equation}
Therefore by the Cauchy theorem it is enough to know the zeros of $g$.
Let us show how these ideas apply in the case of a non-differential operator, a pure differential one and in general when the boundary is a
flat manifold..

Let $a_\gamma(y,D_y,D_\theta)=b(y,D_\theta)$. In this case $\varphi_{\xi,\lambda}(y)=\delta(y-\lambda)$ and $g(\lambda,\xi)=
b(\lambda,\xi)$. Therefore the l.h.s of (\ref{thierry3}) can be computed 
easily: let $\xi_i(y)$ a family of non-degenerate zeros of $b(
y, \xi)$, that is $b(y,\xi_i(y))=0$, $\partial _\xi b(y\xi_i(y))\neq 0$. 
then:
\begin{eqnarray*} 
u(y,\theta) & = & \frac 1 {2\pi i}\int d\lambda\int_{x-i\infty}^{x+i\infty}d\xi\frac 1 {g(\lambda,\xi)}\hat G(\lambda,\xi)e^{\xi\theta}
\varphi_{\xi,\lambda}(y)\\
 &=& \frac 1 {2\pi i}\int d\lambda\int_{x-i\infty}^{x+i\infty}d\xi\frac 1 {b(\lambda,\xi)}\hat G(\lambda,\xi)e^{\xi\theta}
\delta(y-\lambda)\\
&=&  \frac{1}{2\pi i}\int_{x-i\infty}^{x+i\infty}d\xi\frac 1 {b(y,\xi)}\hat G(y,\xi)e^{\xi\theta}\\
&=& \sum_i\frac 1 {\partial_\xi b(y,\xi)|_{\xi=\xi_i(y)}}\hat G(y,\xi_i(y))e^{\xi_i(y)\theta}
\end{eqnarray*}
One sees immediatly on this example that the boundary layer will exist
 as soon as  one of the $\xi_i(y)$ satisfies ${\cal R}e(\xi_i(y)) \ge c >
 0,$  and there is no crossing between the different
branches. 

Let us suppose now that $a_\gamma(y,D_y,D_\theta)=c(D_y,D_\theta)$. In this case $\varphi_{\xi,\lambda}(y)=e^{i\lambda y}$ and $g(\lambda,\xi)=
c(\lambda,\xi)$. Let as before $\xi_i(\lambda)$ a family of 
non-degenerate zeros of $b(\lambda, \xi)$, that is
 $b(\lambda,\xi_i(\lambda))=0$,
 $\partial _\xi b(\lambda,\xi_i(\lambda))\neq 0$. We get: 

\begin{eqnarray*} 
u(y,\theta) & = & \frac 1 {2\pi i}\int d\lambda\int_{x-i\infty}^{x+i\infty}d\xi\frac 1 {g(\lambda,\xi)}\hat G(\lambda,\xi)e^{\xi\theta}
\varphi_{\xi,\lambda}(y)\\
 &=& \frac 1 {2\pi i}\int d\lambda\int_{x-i\infty}^{x+i\infty}d\xi\frac 1 {c(\lambda,\xi)}\hat G(\lambda,\xi)e^{\xi\theta}
e^{i\lambda y}\\
&=& \sum_i\frac 1 {\partial_\xi b(\lambda,\xi)|_{\xi=\xi_i(\lambda)}}\hat G(\lambda,\xi_i(\lambda))e^{\xi_i(\lambda)\theta}
e^{i\lambda y}\\.
\end{eqnarray*}
Therefore the boundary layer will exist as soon as
  there exists a family $\xi_i(\lambda)$ 
with ${\cal R}e(\xi_i(\lambda)) \ge  c > 0,\ \forall \lambda$.

In the case where the boundary is $\R^n$ (or $\T^n$) and that the operator $a_\gamma(y,D_y,D_\theta)$ is a differential operator in $y$ with
polynomial coefficients a decomposition similar to the preceding one can
be handled with the so-called coherent states. 
Let us briefly describe  the construction on $\R$, the generalization
 to $\R^n$ (and $\T^n$ by periodization) is straightforward.
Let $(p,q)\in\R^2$ and let $\varphi_{p,q}$ the family of functions defined
in \cite{Folland:1989}
\[
\varphi_{p,q}(y):=\pi^{-1/4}e^{-(y-q)^2/2}e^{ipy}.
\]
Let us now consider the operator $a_\gamma(y,D_y,\xi)=a(y,D_y,\xi)$. Being
with polynomial coefficients it can be surely rewritten as
\begin{equation}\label{thierry23}
a(y,D_y,\xi)=\sum_{k,j}\alpha_{k,j} (\xi)(y+\partial_y)^k(y-\partial_y)^j
\end{equation}
Let us define $a_{AW}(q,p,\xi) := \sum_{k,j}\alpha_{k,j} (\xi)(q+ip)^k(y-ip)^j$
\begin{lem}
$$a(y,D_y,\xi)=\int_{\R^2} a_{AW}(q,p,\xi) |\varphi_{p,q}><\varphi_{p,q}|dpdq$$
where the integral is understood in the weak sense.
\end{lem}

{\bf Proof}: an easy computation shows that we have the following
decomposition of  identity
\[
\sum_{\R^2}|\varphi_{p,q}><\varphi_{p,q}|dpdq=1
\]
Moreover it is easy to check that
\[
(y+\partial_y)\varphi_{p,q}=(q+ip)\varphi_{p,q}
\]
and therefore
\[
<\varphi_{p,q}|(y-\partial_y)=(q-ip)<\varphi_{p,q}|.
\]
Inserting in (\ref{thierry23}) the decomposition of the identity we get:
\begin{eqnarray*}
(y,D_y,\xi) & = &  \sum_{k,j}\alpha_{k,j}
  (\xi)(y+\partial_y)^k(y-\partial_y)^j\\
& = & \sum_{k,j}\alpha_{k,j} (\xi)(y+\partial_y)^k\int_{\R^2}|\varphi_{p,q}><\varphi_{p,q}|dpdq(y-\partial_y)^j\\
& = & \int_{\R^2}\sum_{k,j}\alpha_{k,j}(q+ip)^k|\varphi_{p,q}><\varphi_{p,q}|(q-ip)^j dpdq\\
& = & \int_{\R^2}a_{AW}(q,p,\xi)|\varphi_{p,q}><\varphi_{p,q}| dpdq. 
\end{eqnarray*} 
It is easy to check that the same argument as for the preceding case holds,
and that the  boundary layer exists as soon as it
exists a family of non-degenerate zeros, $a_{AW}(q,p,\xi(q,p))=0$ with
 real part strictly positive (uniformly in $(q,p)$).
A detailed description of this method will be given
elsewhere.

\subsubsection{Almost solvability}

Contrary to what previous cases suggest, the solvability of
\eqref{equationbl} can not be deduced in general from a symbolic analysis. 
 Indeed, to go from a local analysis
(with a local symbol near $y \in Y'$) to a global one (on all $Y'$)
 requires a ``patching process'', involving  a covering of the boundary and
a partition of unity. Such process introduces regularizing operators,
responsible for a loss of information at low frequencies. 
One can clarify this idea in the light of classical results by
 Tr\`eves \cite[Chapitre 3]{Tre}. These results were obtained to describe
the regularity properties of elliptic operators. However, they give some
insight into our boundary layer problem. Let us
consider the ``tangential principal symbol''
$$ A_{\gamma}(y, D_y, D_\theta) \; := \: \sum_{k=0}^m \sigma(a_k)(y,D_y) \:
D_\theta^k, $$
where $a_\gamma := \sum_{k=0}^m a_k(y, D_y) \, D_\theta^k$.
We suppose that 
\begin{description}
\item[i)] $Y'$ is compact
\item[ii)] There exists $q > 0$, such that 
$$\mbox{ degree }\sigma(a_k) \: = \: \mbox{ degree }\sigma(a_m) \, + \,
q(m-k), \quad \forall k. $$ 
\item[iii)] There are two integers, $m^+$ and $m^-$ such that $m^+ + m^- =
  m$, $m^+ \ge 1$, and for all $(y ,\zeta)$ in $T^* Y' - \{0\}$, the
  polynomial $A_\gamma(y, \zeta, \eta)$ with respect to $\eta$ has exactly
  $m^+$ roots with positive imaginary  part, $m^-$ roots with negative
  imaginary part. 
\end{description}
Under these assumptions, we state 
\begin{theo} \label{regularizing}
There exists a family $R_k$, $k=0, \ldots, m$ of regularizing operators
over $Y'$ such that
$$ \tilde{a}_\gamma(y, D_y, D_\theta) \: := \: a_\gamma(y, D_y, D_\theta) \:
+ \: \sum_{k=0}^m R_k(y, D_y) \, D_\theta^k $$
satisfies:  $ \left\{ u \in {\cal X}, \: \tilde{a}_\gamma u
= 0 \right\} \: \neq \: \{ 0 \}$, and 
for all $f \in {\cal X}$, equation $\tilde{a}_\gamma u = f$ has a
solution in ${\cal X}$. 
\end{theo}
Broadly speaking, this theorem says that a singular operator whose ``symbol
has good properties'' (assumptions i) and ii)) is ``close'' to a boundary layer
operator, {\it i.e.} up to additional regularizing operators. However, as these
regularizing operators do not vanish, one can not conclude to the
solvability of \eqref{equationbl} in ${\cal X}$. 

\medskip
{\bf Proof:} We start noticing that 
$$a_\gamma(y, D_y, D_\theta) \: = \: a_m(y, D_y) \, D_\theta^m \: + \: 
\sum_{k=0}^{m-1} a_k(y, D_y) \, D_\theta^k$$ 
with $a_m$ elliptic. Up to
compose by a parametrix of $a_m$, which would yield an additional $R(y,
D_y) D_\theta^m$ term with $R$ regularizing, we can assume that $a_m = 1$. 

\medskip
This simplification, together with ii) and iii) enter the framework studied
in \cite{Tre}. Tr\`eves considered ii) with $q=1$, but its arguments extend
straightforwardly to $q \in \N^*$.  
In particular, we have the ``almost factorization'' 
(see \cite[pages 157-162]{Tre}):
$$ a_\gamma(y, D_y, D_\theta) \: = \: a_\gamma^-(y, D_y, D_\theta) \: 
   a_\gamma^+(y, D_y, D_\theta) \: + \: R(y, D_y, D_\theta), $$
where 
$\displaystyle 
\:  a_\gamma^\pm := \sum_{k=0}^{m^\pm} a_k^\pm(y, D_y) \, D_\theta^k, \quad 
R := \sum_{k=0}^{m^+} R_k(y, D_y) \, D_\theta^k.$\\
The $a_k^\pm$, resp. $R_k$, are pseudo-differential operators of degree 
$q(m^\pm - k)$ resp. $-\infty$.  Moreover, for all $(y, \zeta)$, the roots of 
$$ A_\gamma^\pm(y, \zeta, \eta) \: := \: \sum_{k=0}^{m^\pm} 
\sigma(a_k^\pm)(y, \zeta) \, \eta^k $$
are the $m^\pm$ roots of $A_\gamma(y, \zeta, \cdot)$ 
with $\gtrless$ imaginary part. Thus,
from above decomposition, we can rewrite equation \eqref{equationbl} as a
triangular system:
\begin{equation} \label{triangular}
\begin{aligned}
& a_\gamma^+ v \: + \: R u \: = \: f, \\
& a_\gamma^- u \: = \: v 
\end{aligned}
\end{equation}
Still following \cite{Tre}, we denote 
$-\Delta_y$ the Laplace-Beltrami operator on $Y'$ and set
$\Lambda = (1 - \Delta_y)^{1/2}$. Through the change of variables 
 \begin{equation} 
\begin{aligned}
&D_\theta v^j \: = \: \Lambda^q v^{j+1}, \quad j = 1, \ldots, m^+ - 1, \\
&D_\theta u^j \: = \: \Lambda^q u^{j+1}, \quad j = 1, \ldots, m^- - 1, 
\end{aligned}
\end{equation} 
we can write \eqref{triangular} as a first order system 
\begin{equation} 
\begin{aligned}
& D_\theta V \,  \: - \: {\cal A}^+_\gamma V \: = \: {\cal R} U \: + \: F, \\
& D_\theta U  \, \: - \: {\cal A}^-_\gamma U \: = \: {\cal J} \, V 
\end{aligned}
\end{equation}
with $V = (v^j)$, $\: U = (u^j)$, 
$${\cal R} U = (0, \ldots, \Lambda^{1-m^+} R
u^1), \quad F = (0, \ldots, \Lambda^{1-m^+} f), \quad{\cal J} V = (0, \ldots,
\Lambda^{1-m^-} v^1).$$ The $ {\cal A}_\gamma^\pm$ are matricial
pseudo-differential operators of order $q$ over $Y'$. Moreover, they satisfy
$$ \det\left( \eta - \sigma({\cal A}_\gamma^\pm)(y, \zeta) \right) \: = \: 
A_\gamma^\pm(y, \zeta, \eta).$$

\medskip
Hence, we can apply \cite[theorem 1.1, page 134]{Tre} on parametrices for
parabolic type operators: one can find a function ${\cal U}^\pm(\theta)$,
with values in matricial pseudo-differential operators of order zero over
$Y'$, such that 
$$ D_\theta \, {\cal U}^\pm \: \mp {\cal A}_\gamma^\pm \, {\cal U}^\pm \: 
\sim \: 0 \mbox{ in } Y' \times \{ \theta > 0 \}, \quad 
{\cal U}^\pm\vert_{\theta= 0} = I^\pm $$
where $I^\pm$ is the identity of $\C^{m^\pm}$. 
In each local chart $\displaystyle 
\left({\cal O}', \chi'\right)$ of $Y'$, the local
symbol $\displaystyle 
(y, \zeta) \: \mapsto \:  {\cal U}^\pm_{\chi'}(y, \theta, \zeta) \:$ of ${\cal
U}^\pm(\theta)$ satisfies 
\begin{description}
\item[i)] ${\cal U}^\pm_{\chi'}$ is $\displaystyle
{\cal C}^\infty$ on $\displaystyle  {\cal O}' \times \R_+ \times \R^{n-1}$. 
\item[ii)] To every compact subset ${\cal K}$ of ${\cal O}'$, to every pair
of $(n-1)$-tuples $\alpha, \beta \in \Z^{n-1}_+$, and to every pair of integers
$r, N \ge 0$, there is a constant $C > 0$ such that for all $y$ in
${\cal K}$, $\theta > 0$,  $\zeta \in \R^{n-1}$, 
\begin{equation} \label{schwartz}
 \| \pa^\alpha_y \, \pa^\beta_\zeta \, \pa^r_\theta \,  {\cal
U}^\pm_{\chi'}(y, \theta, \zeta) \| \: \le \: C \, \theta^{-N} \, (1 +
|\zeta|)^{rq - |\beta| - N q}
\end{equation}
\end{description}
We quote that in theorem 1.1 of Tr\`eves, estimate \eqref{schwartz} holds
only locally in $\theta$, because the symbol ${\cal A} = {\cal A}(\theta)$
under consideration has only local control in $\theta$. In our case, ${\cal
A}^\pm_\gamma$ are independent on $\theta$, so that \eqref{schwartz}
extends to $\theta$ in $\R_+$. 

\medskip
Suppose now that $f \in {\cal X}$ so that $F \in {\cal X}^{m^+}$. 
Let $V^0 \in \left( {\cal C}^\infty(Y')\right)^{m_+}$. Set 
\begin{equation} 
\begin{aligned}
& V \: = \: {\cal U}^+(\theta) \, V^0 \: - \: \int_0^\theta
 {\cal U}^+(\theta -s) \, F(s) \, ds, \\
& U \: = \: - \int_\theta^{+\infty} {\cal U}^-(s-\theta) \, {\cal J} V(s)
\, ds. 
\end{aligned}
\end{equation}
By estimates \eqref{schwartz}, it is clear that $(U, V) \in {\cal X}^m$,
and as ${\cal U}^\pm$ are parametrices, 
\begin{equation} 
\begin{aligned}
& D_\theta V  \: - \: {\cal A}^+_\gamma V \: = \: {\cal R} U \: + \: 
{\cal R}^+ \, V  \: + \:  F, \\
& D_\theta U  \: - \: {\cal A}^-_\gamma U \: = \: {\cal J} \, V \: + \:
{\cal R}^- U 
\end{aligned}
\end{equation}   
for regularizing operators ${\cal R}^\pm$. Back to original variables, we
get the result. $\Box$

\section{Boundary layer expansions} \label{blexpansions}
 \subsection{W.K.B. Ansatz}

The notions and hypothesis introduced above  are linked to 
 boundary layer expansions for
solutions of \eqref{elliptique}-\eqref{elliptiquebc}. Let $Y^1$ to $Y^N$
the connected components of $Y = \pa \Omega$. 
To each $i=1, \ldots,N$, we can associate
$r_i \in \N$ boundary layer exponents $\gamma^i_1 < \ldots <
\gamma^i_{r_i}$. Remind that these exponents are positive rational
numbers. Hence, replacing for simplicity  $\eps$ by $\eps^\delta$ with some
appropriate $\delta > 0$, we can assume that $\delta =1$, 
 $\gamma^i_j \in \N^*, \:
\forall i,j$.

\medskip
As usual in WKB approach, we inject approximations
\eqref{ansatz1}-\eqref{ansatz3} into \eqref{elliptique}. The resulting
equations are expanded and ordered according to powers of $\eps$, and
coefficients of the different powers of $\eps$ are set equal to zero. It
leads to a family of equations on the $u^k$ and $v^{i,k}_j$. To lighten
notations, we set $u^k = 0$,  $v^{i,j}_k = 0$ for $k < 0$. 
\begin{description}
\item[i)] {\em Interior terms}. Denoting $a^\eps(x, D_x) := \sum_{l=0}^h
\eps^l a_l(x, D_x)$, we have  
\begin{equation} \label{profileint}
  \: \sum_{l=0}^h \, a_l(x, D_x) \, u^{k-l} \: = \:  \delta_{0k}\, f,
 \quad \forall k \in \N. 
\end{equation}
 where $\delta_{0k}$ is the Kronecker symbol. This reads in short 
\begin{equation} \label{limitsystem} 
 {\cal A}(x, D_x) U = F 
\end{equation}
where $U := (u^k) \in {\cal C}^\infty(\overline{\Omega})^\N$, $F := (f, 0,
\ldots)$, and 
$${\cal A} : {\cal C}^\infty(\overline{\Omega})^\N \mapsto 
{\cal C}^\infty(\overline{\Omega})^\N,  \quad 
({\cal A}U)^k := \sum a_l u^{k-l}.$$  
\item[ii)] {\em Boundary layer terms}. We deduce from \eqref{localsymbol3} 
\begin{equation} \label{profilebl}
 D^{n^i_j}_\theta \, a_{\gamma^i_j}(y, D_y, D_\theta) \, v^{i,k}_j \: +
\: L(y, \theta, D_y, D_\theta)\left( v^{i,0}_j, \ldots, v^{i,k-1}_j \right)
= 0 
\end{equation}
for $y \in Y^i$ and $\theta \in \R_+$. $L$ is the linear differential operator
coming from lower order terms, and $n^i_j$ is an integer. $L$ has smooth
coefficients, polynomial
in $\theta$. As $D_\theta$ is invertible from ${\cal X}^i$ to 
${\cal X}^i$, we rewrite such equations as:
\begin{equation} \label{profilebl1} a_{\gamma^i_j}(y, D_y, D_\theta) \, v^{i,k}_j \: +
\: M(y, \theta, D_y, D_\theta)\left( v^{i,0}_j, \ldots, v^{i,k-1}_j \right)
= 0. 
\end{equation}
In particular, 
\begin{equation} \label{profilebl2}
a_{\gamma^i_j}(y, D_y, D_\theta) \, v^{i,0}_j \: = \: 0 
\end{equation}
Moreover, as $\gamma^i_j$ is a boundary layer exponent, equation
\eqref{profilebl2} has a non-trivial solution in $\left({\cal
X}^i\right)$. Recursively, equation \eqref{profilebl1} (hence equation
\eqref{profilebl}) has a solution in  $\left({\cal X}^i\right)$. Again,
equations \eqref{profilebl} can be written
$$ {\cal A}_{bl}(y, \theta, D_y, D_\theta) V \: = \: 0, 
\quad V = \left(v^{1,k}_j, \ldots, v^{N,k}_j \right)_{j,k}
 \in \left( \prod_{i=1}^N  {\cal X}^i\right)^\N. $$

\end{description}

\begin{rem}
No other scalings in boundary layer expansion \eqref{ansatz3} 
could lead to ``reasonable'' localized solutions.  This means
that, in view of factorization \eqref{factor2}, any boundary layer term of
the type $v(y, t/\eps^\gamma)$ where $\gamma$ is not a singular exponent
yields equations of the type  
$$ e(y, D_y) \, D_\theta^\alpha v \: = \: f $$
for some elliptic operator $e(y, D_y)$ and $\alpha \in \N$. Except for some
spurious source terms in the case where $e(y, D_y)$ is non-invertible, such
equation has no nontrivial localized solution.
\end{rem}   

\subsection{Assumption (H5) and proof of Theorem \ref{theo1}}
At this point, the boundary layer and regular parts have been described
independently from one another. They are of course linked through the
boundary conditions  \eqref{elliptiquebc}. It can be written 
\begin{equation} \label{elliptiquebcbis}
\sum_{k=0}^{m'} \, B^\eps_{kl}(y, D_y) \, D_\theta^k u\vert_{t=0} \: = \: g_l,
\quad y \in Y, 
\end{equation}
where $ B^\eps_{kl}(y, D_y)$ are smooth differential operators, polynomial
in $\eps$. Let us introduce
$$ \partial {\cal V} \: = \:  \left\{ \left( D^k_\theta V\vert_{t=0}
\right), \quad k=0, \ldots, m', \quad {\cal A}_{bl} V = 0 \right\} \:
\subset \: \left( \prod_{i=1}^N  {\cal C}^\infty(Y^i)^\N
\right)^{m'+1} $$ 
Again, we can inject the Ansatz \eqref{ansatz1}-\eqref{ansatz3} into
\eqref{elliptiquebc}.  This can be written as 
$$ \sum_{k=0}^m {\cal B}_{kl}(y, D_y) D_\theta^k U\vert_{t=0} \: + \: {\cal
  C}_l(y, D_y) W \: = \: G_l, $$
for some differential operators 
${\cal B}_{kl}: {\cal C}^\infty(Y)^\N
\mapsto {\cal C}^\infty(Y)^\N$, $W \in \partial {\cal V}$, and 
$${\cal C}_l: \partial {\cal V} \mapsto  \left( \prod_{i=1}^N {\cal C}^\infty(Y^i)^\N
\right)^{m+1}.$$

\medskip
Introducing any projector $\Pi_l$ on the range of ${\cal C}_l$, we can
rewrite previous equation as
\begin{equation} \label{limitbc} 
\sum_{k=0}^m \Pi_l \, {\cal B}_{kl}(y, D_y) D_\theta^k U\vert_{t=0}  \: =
\: \Pi_l G_l
\end{equation}
The existence of the regular part resumes to:

\mspace
{\bf (H5) One can find $U \in {\cal C}^\infty(\overline{\Omega})^\N$ such
  that 
\begin{equation} \label{limit}
\left\{
\begin{aligned}
& {\cal A}(x, D_x) U = F, \\
& \sum_{k=0}^m \Pi_l \, {\cal B}_{kl}(y, D_y) D_\theta^k U\vert_{t=0}  \: =
\: \Pi_l G_l
\end{aligned}
\right.
\end{equation} } 

\medskip
As assumptions (H1) to (H5) are fulfilled, one can find $u^\eps$ of type
  \eqref{ansatz1}-\eqref{ansatz3} satisfying 
\begin{equation} 
\left\{
\begin{aligned}
& a^{\eps}(x, D_x) \, u^\eps \: \sim \: f, \quad x \in \Omega, \\
&\, b^\eps_{l}(x, D_x) \, u \: \sim \: g_l,
\quad x \in Y, \quad l \in \{1, \dots, L\},
\end{aligned}
\right.
\end{equation}
uniformly on the compact subsets. Theorem \ref{theo1} follows. 

\medskip
One can see system \eqref{limit} as the limit system of
\eqref{elliptique}-\eqref{elliptiquebc}. It means that the small-scale
effects due to boundary layers are contained in this large-scale system. 
 A strong difficulty is that it is {\em infinite-dimensional}:
correction terms $u^i$ at any order may interact. Hence, in most
situations, one does not  treat directly system \eqref{limit}. The
classical attempt is to turn it
 into  {\em an infinite collection of finite-dimensional 
systems}, that are of the same type and solved recursively.This will be
 illustrated on the quasigeostrophic equation  in the next section.  
 
\section{Application to the quasigeostrophic model} \label{sectionqg}
The stationary linearized quasigeostrophic equation reads, for $x$ in a
two-dimensional domain $\Omega$ 
\begin{equation} \label{qg}
\beta \, \pa_{x_1} \psi \: + \: \na^\bot \psi \cdot \na \eta_B \: + \:  r
\Delta \psi \: - \: \frac{1}{{\rm Re}} \Delta^2 \psi \: = \:  \beta \, \curl \tau
\end{equation}
where $\psi  =\psi(x) \in \R$, $x = (x_1, x_2)$ in cartesian coordinates,
$\beta$, $r$ and ${\rm Re}$ are positive parameters, and $\eta_B =
\eta_B(x)$, $\tau = \tau(x)$ are smooth functions.
 It is a celebrated basic model for oceanic circulation: $\psi$ is a stream
 function, associated to horizontal velocity field $u = \na^\bot \psi$,
 $\eta_B$ is a bottom topography term. The expression $\beta \, \pa_{x_1}
 \psi$ comes from the variation of the Coriolis force with latitude, $r
 \Delta \psi$ and $- 1/{\rm Re} \,  \Delta^2 \psi$ are dissipative terms due
 resp.  to friction and viscosity. Finally, $\beta \, \curl \tau$ is a
 vorticity term created by the wind. We refer to \cite{Ped} for all
 necessary details. 
Equation \eqref{qg} is fulfilled with Dirichlet boundary conditions
\begin{equation} \label{qgbc}
\psi = \pa_n \psi = 0, \quad x \in \pa \Omega. 
\end{equation}  
We will apply previous analysis to system \eqref{qg}, \eqref{qgbc}, for
several domains $\Omega$ and ranges of parameters $\beta$, $r$ and ${\rm
  Re}$.  

\subsection{Munk layers} \label{munk}
We first investigate the case where 
\begin{equation} \label{scaling1}
\beta \rightarrow +\infty, \quad r, {\rm Re}, \tau \mbox{ given}
\end{equation}
which corresponds to strong forcing by the wind. We consider a domain of
the type 
$$ \Omega \: := \: \left\{ \chi_w(x_2) \le x_1 \le \chi_e(x_2) \right\}, $$
where $\chi_w$ and $\chi_e$ are smooth, with all derivatives bounded. They
describe the western and eastern coasts. We assume that $x_2 \in \R$ or
$x_2 \in \T$. The case of a closed basin will be evoked in last
subsection. 

\medskip
This boundary layer problem enters the framework given in section
\ref{singularperturbation}, with $\eps = \beta^{-1}$. We start with the
derivation of singular exponents for both connected components of $Y := \pa
\Omega$,   
$$ Y_w := \left\{ (\chi_w(x_2), x_2), \quad x_2 \in \R \mbox{ or } \T
\right\}, \:  Y_e := \left\{ (\chi_e(x_2), x_2), \quad x_2 \in \R \mbox{ or } \T
\right\}. $$
{\em Western boundary layers}. For any $y \in Y_w$, global coordinates 
$x'_1 = x_2$, $x'_2 = x_1 - \chi_w(x_2)$ define a local chart near
$y$. In these coordinates, the corresponding symbol satisfies  as
$\displaystyle |\xi_1| + |\xi_2| \rightarrow +\infty$, 
\begin{multline*}
 a_w(x', \xi_1, \xi_2) \: := \: i \,  \eps^{-1} \, \xi_2 \: + \: O(|\xi_1|
 \,  + \,
 |\xi_2| ) 
 \: -  \: r \, \Bigl( \xi_1^2 \, + \, (1 + \chi'^{2}_w) \, \xi_2^2 \\ + \:
 O(|\xi_1| \, + \, |\xi_2|)\Bigr) 
 \: - \: \frac{1}{{\rm Re}} \, \Bigl( \xi_1^4 \, + \, (1 + \chi'^{2}_w)^2
 \xi_2^2 \, + \, O((|\xi_1| \, + \, |\xi_2|)^3)\Bigr), 
\end{multline*}
  with 
$\displaystyle \chi'_w :=
  \chi_w'(x'_1)$. For all $\displaystyle \xi_1 \in \R$, the roots of
  $\displaystyle  a_w(x',\xi_1, \cdot)$ with  $|\xi_2| \rightarrow  +\infty$ satisfy
    $$ i \, \eps^{-1} \, \xi_2 \: \sim \: \frac{1}{{\rm Re}} \, (1 + \chi'^2_w)^2
  \xi_2^4, \quad  \mbox{{\it i.e.} } \:  \xi_2 \: \sim \: \left( \frac{i{\rm Re} }{(1
    + \chi'^2_w)^2}\right)^{1/3} \, \eps^{-1/3}. $$
{\em Eastern boundary layers}. The derivation is similar, with
$x'_1 = x_2$, $x'_2 = \chi_e(x_2) - x_1$. Roots of $a_e(x', \xi_1, \cdot)$
with singular behaviour satisfy
$$ \xi_1 \: \sim \: \left( \frac{- i {\rm Re}}{ (1 +
  \chi'^2_e)^2}\right)^{1/3} \, \eps^{-1/3}. $$   
 Hence, assumptions
(H1) to (H3) are fulfilled, with the singular exponent $\gamma = 1/3$.   

\medskip
To identify singular operators and boundary layer sizes is easy using the
global coordinates. Indeed,
\begin{align*}
&  a_w(x'_1, 0, D_{x'_1}, \frac{D_\theta}{\eps^{1/3}}) \: \sim \:
\eps^{-4/3} \, \left( D_\theta - \frac{1}{{\rm Re}}\left( 1 +
\chi'^2_w\right)^2 \, D_\theta^4 \right),   \\
& a_e(x'_1, 0, D_{x'_1}, \frac{D_\theta}{\eps^{1/3}}) \: \sim \:
\eps^{-4/3} \, \left( -D_\theta - \frac{1}{{\rm Re}}\left( 1 +
\chi'^2_e\right)^2 \, D_\theta^4 \right)
\end{align*}
Thus, the singular operators are 
$$ a_{w,\gamma} :=  D_\theta - \frac{1}{{\rm Re}}\left( 1 +
\chi'^2_w\right)^2 \, D_\theta^4, \quad a_{e,\gamma} :=  -D_\theta -
\frac{1}{{\rm Re}} \left( 1 + \chi'^2_e\right)^2 \, D_\theta^4.   $$
They have order zero coefficients. Denoting $\displaystyle \alpha_{w,e} :=
{\rm Re}^{1/3}/ (1 + \chi'^2_{w,e})^2$,  the roots of their symbols are $\alpha_w
e^{i \pi /6}$, $\alpha_w e^{5 i \pi/6}$,   $\alpha_w e^{3 i \pi/2}$, resp. 
$\alpha_e e^{-i \pi /6}$, $\alpha_e e^{ i \pi/2}$,   $\alpha_e e^{7 i
  \pi/6}$. Proposition
\ref{solv} applies (two roots, resp. one root with positive imaginary
part). Thus, $\eps^{1/3}$ is a boundary layer size on $Y_{w,e}$. 

\medskip
It remains to show the existence of solutions similar to \eqref{ansatz1},
with 
\begin{align*}
  u^\eps_{int}(x) \: \sim & \: \psi^0(x) \: + \: \eps^{1/3} \, \psi^1(x)
\: + \ldots \\
 v^\eps_{bl}(x) \: \sim & \:  \psi^0_w\left( \frac{x_1 -
  \chi_w(x_2)}{\eps^{1/3}}, x_2 \right) \: + \:  \psi^0_e\left( 
\frac{\chi_e(x_2) - x_1}{\eps^{1/3}}, x_2 \right)  \\
& +  \:  \eps^{1/3}  \psi^1_w\left( \frac{x_1 - \chi_w(x_2)}{\eps^{1/3}}, x_2 \right) \: + \:
  \ldots 
\end{align*}
Equation \eqref{qg} yields 
\begin{equation} \label{sverdrup}
\pa_{x_1} \psi^i \: = \: f^i \: + \:  \delta_{0i} \, \curl \tau, 
\end{equation}
where $f^i$ involves $u^k$, $k \le i-1$. Equation \eqref{qgbc} yields 
\begin{align*}
& \psi^i_w\vert_{\theta = 0} \: + \: \psi^i\vert_{Y_w} = 0, \quad 
 \psi^i_e\vert_{\theta = 0} \: + \: \psi^i\vert_{Y_e} = 0, \\
& D_\theta \psi^i_w\vert_{\theta = 0} = g^i_w, \quad  D_\theta
 \psi^i_e\vert_{\theta = 0} = g^i_e, 
\end{align*}
where $g^i_{w,e}$ depends on $u^k$ and $u^k_{w,e}$, $k \le i-1$. The
western boundary layer operator has two characteristic roots with positive
imaginary part. One can find localized solutions $\psi_w$ of $a_{w, \gamma}
\psi_w = f_w$ for arbitrary values of $\psi_w\vert_{\theta=0}$, $D_\theta
\psi_w\vert_{\theta=0}$. This is not the case for the eastern boundary
layer, for which only one condition can be given. The last boundary
condition prescribes $D_\theta \psi^0_e\vert_{\theta=0} = 0$. Recursively,
we easily deduce that  the approximate solutions exist if  
\begin{equation} \label{sverdrup2}
\pa_x \psi = f, \quad \psi\vert_{Y_e} = g, \quad f \in {\cal
  C}^\infty(\overline{\Omega}), \quad g \in {\cal C}^\infty(Y_e) 
\end{equation} 
has a solution. As $\psi(x_1, x_2) = \int_{\chi_e(x_2)}^{x_1} f \: + \:
g(x_1)$ is the unique solution of \eqref{sverdrup2}, this proves the
existence of boundary layer expansions. These boundary layers
are called {\em Munk layers}. 

\subsection{Stommel and friction layers}
We now investigate the case where 
\begin{equation} \label{scaling2}
\beta \rightarrow +\infty, \quad \frac{\beta^{2/3}}{r} \rightarrow 0, \quad
      {\rm Re}, \tau \mbox{ given}, 
\end{equation}
which emphasizes the role of friction. In terms of $\eps := \beta^{-1}$
and $\eps_s := r/\beta$, the last condition reads $\eps^{1/3} / \eps_s
\rightarrow 0$. To fit exactly the context under consideration, we assume that
$\eps_s = \eps^\delta$ for some rational number $\delta < 1/3$. The domain
$\Omega$  remains the same as in previous subsection. 

\medskip
Remind that the western symbol reads 
  \begin{multline*}
 a_w(x', \xi_1, \xi_2) \: := \: i \,  \eps^{-1} \, \xi_2 \: + \: O(|\xi_1|
 \,  + \,
 |\xi_2| ) 
 \: -  \: \eps \, \eps^{-1}  \, \Bigl( \xi_1^2 \, + \, (1 + \chi'^{2}_w) \, \xi_2^2 \\ + \:
 O(|\xi_1| \, + \, |\xi_2|)\Bigr) 
 \: - \: \frac{1}{{\rm Re}} \, \Bigl( \xi_1^4 \, + \, (1 + \chi'^{2}_w)^2
 \xi_2^4 \, + \, O((|\xi_1| \, + \, |\xi_2|)^3)\Bigr).
\end{multline*}
For all $\xi_1 \in \R$, the singular roots of $a_w(x', \xi_1, \cdot)$
satisfy either 
\begin{align*}  i \, \eps^{-1} \, \xi_2 \: \sim \: \eps_s \, \eps^{-1} \, \left( 1 +
\chi'^2_w\right) \xi_2^2, \: & \mbox{ {\it i.e.} } \, \xi_2 \: \sim \: 
\frac{i}{1 + \chi'^2_w} \, \eps_s^{-1}, \quad \mbox{ or} \\ 
 \frac{\eps_s}{\eps} \left( 1 +
\chi'^2_w\right) \xi_2^2 \: \sim \: \frac{-1}{{\rm Re}}  \left( 1 +
\chi'^2_w\right)^2   \xi_2^4, \:  & \mbox{ {\it i.e.} } \, \xi_2 \: \sim \: 
\frac{(-{\rm Re})^{1/2}}{ 1 + \chi'^2_w} \,  \left( \frac{\eps_s}{\eps}
\right)^{1/2}.  
\end{align*}
Similarly, for the eastern symbol, 
$$ \xi_2 \: \sim \: \frac{-i}{1 + \chi'^2_w} \, \eps_s^{-1}, \quad \mbox{
  or} \: \xi_2 \: \sim \: 
\frac{(-{\rm Re})^{1/2}}{ 1 + \chi'^2_e} \,  \left( \frac{\eps_s}{\eps}
\right)^{1/2}.
$$
There are two singular exponents, uniformly on $Y$: $\gamma_1 = \delta$
({\it i.e.} $\eps^{\gamma_1} = \eps_s$) and $\gamma_2 = \frac{1-\delta}{2}$ 
({\it i.e.} $\eps^{\gamma_2} = \sqrt{\eps/ \eps_s}$). 

\medskip
The singular operators read 
\begin{align*}
& a_{w, \gamma_1} \, = \, i \, D_\theta \, - \, r D_\theta^2, \quad 
a_{e, \gamma_1} \, = \, -i \, D_\theta \, - \, r D_\theta^2, \\
& a_{w, \gamma_2} \, = \, -r Id \, -\, \frac{1}{{\rm Re}} D_\theta^2, \quad 
  a_{e, \gamma_2} \, = \, -r Id \, -\, \frac{1}{{\rm Re}} D_\theta^2.
\end{align*}
Thanks to proposition \ref{solv}, we deduce that $\eps^{\gamma_2}$ is a
boundary layer size on $Y_{w,e}$, and  $\eps^{\gamma_1}$ is only a boundary
layer size on $Y_e$ (one root with positive imaginary part for
$a_{w,e}^{\gamma_2}$ and $a_w^{\gamma_1}$, no such root for
  $a_e^{\gamma_1}$). 

\medskip
We look for approximate solutions with 
\begin{align*}
 u^\eps_{int} \: &  = \:  \psi^0  \: + \: \eps^{\gamma_2 - \gamma_1} \,
\psi^1 \: + \: \eps^{2(\gamma_2 - \gamma_1)} \, \psi^2 \: + \: \ldots, \\
 v^\eps_{bl} \: & = \: \psi^0_{w,1}\left(\frac{x_1 -
  \chi_w(x_2)}{\eps^{\gamma_1}}, x_2 \right) \: + \:   \psi^0_{w,2}\left(\frac{x_1 -
  \chi_w(x_2)}{\eps^{\gamma_2}}, x_2 \right) \\
& + \: \psi^0_{e,2}\left(\frac{
  \chi_e(x_2)-x_1}{\eps^{\gamma_2}}, x_2 \right) \: + \: \eps^{\gamma_2 -
   \gamma_1} \psi^1_{w,1}\left(\frac{x_1 -
  \chi_w(x_2)}{\eps^{\gamma_1}}, x_2 \right) \: + \: \ldots
 \end{align*}
Equation \eqref{qg} yields \eqref{sverdrup}. Equation \eqref{qgbc} yields
\begin{align*}
& \psi^i_{w,1}\vert_{\theta = 0} \: + \:  \psi^i_{w,2}\vert_{\theta = 0} 
\: + \: \psi^i\vert_{Y_w} = 0, \quad 
 \psi^i_{e,2}\vert_{\theta = 0} \: + \: \psi^i\vert_{Y_e} = 0, \\
& D_\theta \psi^i_{w,2}\vert_{\theta = 0} = g^i_w, \quad  D_\theta
 \psi^i_{e,2}\vert_{\theta = 0} = g^i_e, 
\end{align*}
where $g^i_{w,e}$ involve $u^k$ and $u^k_{w,e,1,2}$ for $k \le i-1$. All
boundary layer operators have one root with positive imaginary part, thus
only one condition can be given at the boundary. The last line prescribes
$D_\theta \psi^0_{w,e,2} = 0$. Reasoning recursively, the existence of 
 approximate solutions  relies on the solvability of
 \eqref{sverdrup2}. The boundary layers with size $\eps_s$ and $\eps_f :=
 \sqrt{\eps_s/\eps}$  are known as the Stommel and friction layers. 
\subsection{Geostrophic degeneracy}
We consider again the asymptotics \eqref{scaling1}. In subsection
\ref{munk}, we have considered a domain with western and eastern
boundaries, but no northern or southern boundaries. Let us examine here the
case of a closed basin. We assume simply that 
$$ \Omega  \: = \: \left\{ x_1^2 + x_2^2 \: \le 1  \right\},
 \quad Y = \pa \Omega = \mathbb{S}^1. $$
Let us look for singular exponents near $y_0 := (0,1)$. We chose local
coordinates $x'_1 = x_1$, $\displaystyle 
x'_2 = \sqrt{1 - x'^2_1} -  x_2$. The local
symbol reads, as $|\xi_1| + |\xi_2| \rightarrow \infty$  
\begin{align*}
a^\eps_{\chi}(x', \xi_1, \xi_2) \: = &
\: i \, \eps^{-1}  \, \left( \xi_1  -
\frac{2 x'_1}{\sqrt{1 - x'^2_1}} \, \xi_2 \right)  \\
& \:  + \: 
i \left( \begin{smallmatrix} \xi_1 - \frac{2 x'_1}{\sqrt{1 - x'^2_1}} \xi_2
  \\ -\xi_2 \end{smallmatrix} \right)^{\bot} \cdot 
 \left( \begin{smallmatrix} \pa_{x'_1} \eta_B  - 
\frac{2 x'_1}{\sqrt{1 - x'^2_1}} \pa_{x'_2} \eta_B
  \\ - \pa_{x'_2} \eta_B \end{smallmatrix} \right)  \\
& \:  -\:  
r \, \left( \xi_1^2 + \left( \frac{4 x'^2_1}{1 - x'^2_1} + 1 \right) \,
\xi_2^2 \: + \: O(|\xi_1| + |\xi_2|)\right)  \\
& \:  -\: 
\frac{1}{{\rm Re}} \left( \xi_1^4 \, + \,  \left( \frac{4 x'^2_1}{1 -
   x'^2_1} + 1 \right)^2 \xi_2^4 \right) \, + \, O((|\xi_1| + |\xi_2|)^3) 
\end{align*}
Let $\xi_1 \in \R^*$. For $y \neq y_0$, {\it i.e.} $x'_1 \neq 0$, the 
roots $\xi_2$ with $|\xi_2| \rightarrow +\infty$ satisfy 
$$  2 i \, \eps^{-1}  \, 
 \frac{1}{{\rm Re}}  \: \frac{x'_1}{\sqrt{1 - x'^2_1}} \, \xi_2  \: \sim \: 
\left( \frac{4 x'^2_1}{1 - x'^2_1} + 1 \right)^2 \xi_2^4, $$
which yields $\gamma = 1/3$ in agreement  with subsection \ref{munk}. For $y
= y_0$, {\it i.e.} $x'_1 = 0$, this relation degenerates. The roots $\xi_2$
going to infinity satisfy 
$$ i \, \eps^{-1} \, \xi_1 \: \sim \: \frac{1}{{\rm Re}} \, \xi_2^4, \:
{\it i.e.} \: \xi_2 \: \sim \: \left( i \, \xi_1 {\rm Re} \right)^{1/4} \,
\eps^{-1/4} $$ 
so that $\gamma = 1/4$. In particular, there are no singular exponents
uniformly on a neighborhood of  $y_0 = (0,1)$ and (H3) fails. This
phenomenon is known as geostrophic degeneracy. It takes place in various
systems of fluid mechanics, including rotating fluids (as can be seen on
figure \ref{figure}, near the equator of the inner sphere,
 following \cite{Ste2}). To understand the
structure of the solutions near such ``turning points'' is an open
question, both from physical and mathematical viewpoints.
\bibliographystyle{plain}
\def\cprime{$'$}

\end{document}

%% file: sphere.eepic
\setlength{\unitlength}{0.00087489in}
\begingroup\makeatletter\ifx\SetFigFont\undefined%
\gdef\SetFigFont#1#2#3#4#5{%
  \reset@font\fontsize{#1}{#2pt}%
  \fontfamily{#3}\fontseries{#4}\fontshape{#5}%
  \selectfont}%
\fi\endgroup%
{\renewcommand{\dashlinestretch}{30}
\begin{picture}(5659,7869)(0,-10)
\put(405.000,3882.000){\arc{3600.000}{4.7124}{7.8540}}
\path(2475,7842)(2475,4332)
\path(2475,3477)(2475,57)
\path(2025,7842)(2025,5232)
\path(1935,2442)(1935,102)
\path(2475,6492)(2475,7842)
\path(3105,7842)(3105,147)
\path(3825,7842)(3825,282)
\texture{88555555 55000000 555555 55000000 555555 55000000 555555 55000000 
	555555 55000000 555555 55000000 555555 55000000 555555 55000000 
	555555 55000000 555555 55000000 555555 55000000 555555 55000000 
	555555 55000000 555555 55000000 555555 55000000 555555 55000000 }
\dashline{60.000}(2475,4287)(1305,4287)
\dashline{60.000}(2430,3432)(1305,3432)
\path(1380.000,4077.000)(1350.000,4197.000)(1320.000,4077.000)
\path(1350,4197)(1350,3477)
\path(1350,4197)(1350,3477)
\path(1320.000,3597.000)(1350.000,3477.000)(1380.000,3597.000)
\path(2085.000,3852.000)(2205.000,3882.000)(2085.000,3912.000)
\path(2205,3882)(1755,3882)
\path(2205,3882)(1755,3882)
\path(2865.000,3912.000)(2745.000,3882.000)(2865.000,3852.000)
\path(2745,3882)(3060,3882)
\path(2745,3882)(3060,3882)
\path(2577.442,3767.678)(2475.000,3837.000)(2532.768,3727.625)
\path(2475,3837)(4815,1227)
\path(2475,3837)(4815,1227)
\path(930.000,1962.000)(900.000,2082.000)(870.000,1962.000)
\path(900,2082)(900,867)
\path(900,2082)(900,867)
\path(1935,237)(1935,12)
\path(1935,237)(1935,12)
\path(3105,642)(3105,12)
\path(3105,642)(3105,12)
\path(2475,507)(2475,12)
\path(2475,507)(2475,12)
\path(3825,597)(3825,12)
\path(3825,597)(3825,12)
\path(2145.000,5982.000)(2025.000,5952.000)(2145.000,5922.000)
\path(2025,5952)(2475,5952)
\path(2025,5952)(2475,5952)
\path(2355.000,5922.000)(2475.000,5952.000)(2355.000,5982.000)
\path(2595.000,6612.000)(2475.000,6582.000)(2595.000,6552.000)
\path(2475,6582)(3105,6582)
\path(2475,6582)(3105,6582)
\path(2985.000,6552.000)(3105.000,6582.000)(2985.000,6612.000)
\path(3225.000,7242.000)(3105.000,7212.000)(3225.000,7182.000)
\path(3105,7212)(3825,7212)
\path(3105,7212)(3825,7212)
\path(3705.000,7182.000)(3825.000,7212.000)(3705.000,7242.000)
\path(360,6537)(945,5727)
\path(360,6537)(945,5727)
\path(850.421,5806.717)(945.000,5727.000)(899.062,5841.846)
\path(405,1812)(408,1812)(415,1813)
	(427,1814)(444,1816)(468,1818)
	(496,1821)(528,1824)(563,1827)
	(599,1831)(634,1834)(668,1838)
	(700,1841)(731,1845)(759,1848)
	(784,1852)(808,1855)(831,1858)
	(851,1861)(871,1865)(889,1868)
	(908,1872)(929,1877)(950,1882)
	(971,1887)(992,1893)(1014,1899)
	(1035,1906)(1056,1913)(1077,1920)
	(1098,1928)(1119,1935)(1139,1943)
	(1159,1951)(1179,1960)(1197,1968)
	(1215,1976)(1233,1984)(1250,1992)
	(1268,1999)(1285,2008)(1302,2016)
	(1321,2025)(1339,2034)(1358,2043)
	(1378,2053)(1399,2063)(1419,2074)
	(1440,2085)(1461,2096)(1481,2108)
	(1502,2119)(1522,2130)(1541,2142)
	(1559,2153)(1578,2164)(1595,2176)
	(1613,2187)(1628,2197)(1643,2208)
	(1659,2219)(1675,2231)(1691,2243)
	(1708,2256)(1725,2269)(1741,2282)
	(1758,2296)(1775,2310)(1792,2324)
	(1809,2339)(1825,2353)(1841,2367)
	(1857,2382)(1872,2395)(1886,2409)
	(1900,2423)(1914,2436)(1928,2449)
	(1942,2464)(1957,2480)(1972,2496)
	(1988,2512)(2003,2529)(2019,2546)
	(2035,2564)(2050,2582)(2066,2601)
	(2081,2619)(2096,2637)(2110,2656)
	(2124,2674)(2137,2691)(2149,2708)
	(2161,2725)(2172,2741)(2183,2757)
	(2193,2773)(2203,2789)(2212,2806)
	(2222,2823)(2232,2841)(2241,2859)
	(2250,2878)(2259,2896)(2268,2916)
	(2277,2935)(2285,2954)(2293,2973)
	(2301,2991)(2308,3010)(2314,3028)
	(2321,3045)(2327,3062)(2333,3079)
	(2338,3097)(2344,3114)(2349,3133)
	(2355,3152)(2361,3173)(2368,3196)
	(2375,3221)(2382,3248)(2389,3276)
	(2397,3306)(2405,3335)(2412,3363)
	(2419,3388)(2424,3407)(2427,3421)
	(2429,3428)(2430,3432)
\path(405,5997)(408,5997)(415,5996)
	(426,5995)(443,5993)(465,5990)
	(491,5987)(520,5983)(550,5980)
	(580,5976)(609,5972)(637,5969)
	(663,5966)(686,5962)(708,5959)
	(729,5956)(748,5953)(767,5951)
	(785,5948)(803,5944)(822,5941)
	(842,5937)(861,5933)(882,5929)
	(902,5925)(923,5920)(944,5915)
	(966,5910)(987,5905)(1008,5899)
	(1029,5894)(1049,5888)(1068,5883)
	(1087,5877)(1105,5871)(1122,5866)
	(1139,5860)(1155,5854)(1171,5849)
	(1187,5843)(1203,5836)(1219,5830)
	(1236,5823)(1253,5815)(1270,5808)
	(1287,5800)(1304,5791)(1321,5783)
	(1338,5774)(1355,5765)(1371,5757)
	(1387,5748)(1402,5739)(1418,5730)
	(1433,5721)(1448,5712)(1462,5703)
	(1478,5693)(1493,5683)(1510,5672)
	(1526,5661)(1543,5650)(1561,5637)
	(1579,5625)(1596,5612)(1614,5600)
	(1632,5587)(1649,5574)(1665,5561)
	(1681,5549)(1697,5537)(1712,5525)
	(1726,5513)(1740,5502)(1754,5491)
	(1767,5479)(1781,5467)(1795,5455)
	(1808,5442)(1822,5429)(1836,5416)
	(1850,5402)(1864,5388)(1878,5375)
	(1891,5361)(1904,5347)(1917,5334)
	(1929,5320)(1940,5307)(1951,5295)
	(1962,5282)(1973,5269)(1984,5255)
	(1996,5241)(2007,5225)(2019,5209)
	(2032,5193)(2044,5176)(2056,5159)
	(2069,5141)(2081,5123)(2094,5105)
	(2106,5088)(2117,5071)(2128,5054)
	(2139,5038)(2150,5023)(2160,5007)
	(2169,4993)(2178,4979)(2188,4965)
	(2198,4950)(2208,4935)(2218,4919)
	(2228,4902)(2239,4885)(2249,4868)
	(2259,4850)(2270,4832)(2279,4814)
	(2289,4797)(2298,4779)(2307,4761)
	(2316,4743)(2324,4725)(2333,4707)
	(2340,4690)(2347,4673)(2354,4655)
	(2361,4635)(2369,4614)(2377,4591)
	(2386,4566)(2395,4539)(2405,4509)
	(2415,4478)(2425,4446)(2436,4413)
	(2446,4382)(2454,4353)(2462,4329)
	(2468,4310)(2472,4298)(2474,4290)(2475,4287)
\path(2475,4287)(2476,4287)(2479,4284)
	(2488,4278)(2501,4269)(2516,4258)
	(2533,4246)(2548,4235)(2562,4225)
	(2574,4215)(2584,4206)(2594,4198)
	(2602,4190)(2603,4189)(2611,4181)
	(2619,4171)(2628,4161)(2636,4151)
	(2645,4139)(2653,4127)(2661,4115)
	(2668,4103)(2675,4091)(2681,4078)
	(2687,4066)(2692,4055)(2693,4054)
	(2698,4042)(2702,4029)(2707,4016)
	(2712,4001)(2716,3986)(2720,3971)
	(2724,3955)(2728,3940)(2731,3925)
	(2733,3910)(2736,3896)(2737,3882)
	(2738,3882)(2739,3868)(2740,3854)
	(2741,3839)(2742,3824)(2742,3809)
	(2742,3793)(2741,3778)(2740,3763)
	(2738,3748)(2736,3735)(2733,3722)
	(2730,3710)(2730,3709)(2726,3698)
	(2722,3686)(2717,3674)(2711,3661)
	(2705,3649)(2698,3638)(2691,3626)
	(2683,3616)(2676,3606)(2669,3597)
	(2662,3589)(2655,3582)(2647,3574)
	(2638,3567)(2629,3559)(2619,3552)
	(2609,3546)(2599,3539)(2588,3533)
	(2578,3527)(2568,3521)(2558,3514)
	(2549,3509)(2539,3503)(2529,3497)
	(2517,3489)(2503,3480)(2488,3470)
	(2471,3459)(2456,3449)(2443,3440)
	(2434,3435)(2431,3432)(2430,3432)
\put(540,3027){\makebox(0,0)[lb]{\smash{{{\SetFigFont{12}{14.4}{\rmdefault}{\mddefault}{\updefault}Inner Sphere}}}}}
\put(540,3792){\makebox(0,0)[lb]{\smash{{{\SetFigFont{12}{14.4}{\rmdefault}{\mddefault}{\updefault}$E^{1/5}$}}}}}
\put(3150,7347){\makebox(0,0)[lb]{\smash{{{\SetFigFont{12}{14.4}{\rmdefault}{\mddefault}{\updefault}$E^{1/4}$}}}}}
\put(2520,6717){\makebox(0,0)[lb]{\smash{{{\SetFigFont{12}{14.4}{\rmdefault}{\mddefault}{\updefault}$E^{1/3}$}}}}}
\put(4860,1047){\makebox(0,0)[lb]{\smash{{{\SetFigFont{12}{14.4}{\rmdefault}{\mddefault}{\updefault}$E^{2/5}$}}}}}
\put(315,732){\makebox(0,0)[lb]{\smash{{{\SetFigFont{12}{14.4}{\rmdefault}{\mddefault}{\updefault}Ekman Layer}}}}}
\put(0,6717){\makebox(0,0)[lb]{\smash{{{\SetFigFont{12}{14.4}{\rmdefault}{\mddefault}{\updefault}$E^{1/2}$}}}}}
\put(2115,6132){\makebox(0,0)[lb]{\smash{{{\SetFigFont{12}{14.4}{\rmdefault}{\mddefault}{\updefault}$E^{2/7}$}}}}}
\end{picture}
}